\documentclass[12pt, reqno, a4paper]{amsart}

\usepackage{ amssymb, amsmath, enumerate, amsfonts, amsthm, mathrsfs, url, bm, mathtools}
\usepackage[all,cmtip]{xy}

\numberwithin{equation}{section}

\usepackage{xcolor}
	
\usepackage[backref]{hyperref}
\hypersetup{
	colorlinks,
    linkcolor={blue!60!black},
    citecolor={blue!60!black},
    urlcolor={red!60!black}
}

\usepackage[margin=1.25in]{geometry}

\RequirePackage{doi}

\usepackage[square,sort,comma,numbers]{natbib}
 \newcommand{\set}[1]{\left\{#1\right\}}
\newcommand{\bigset}[1]{\bigl\{ #1 \bigr\}}

\newcommand{\abs}[1]{\left| #1\right|}
\newcommand{\bigabs}[1]{\bigl| #1 \bigr|}
\newcommand{\Bigabs}[1]{\Bigl| #1 \Bigr|}

\newcommand{\bigbrac}[1]{\bigl( #1 \bigr)}
\newcommand{\Bigbrac}[1]{\Bigl( #1 \Bigr)}

\newcommand{\norm}[1]{\left\| #1\right\|}
\newcommand{\bignorm}[1]{\big\| #1 \big\|}

\newcommand{\ang}[1]{\left\langle#1\right\rangle}

\newcommand{\N}{\mathbb{N}}
\newcommand{\Z}{\mathbb{Z}}

\newcommand{\R}{\mathbb{R}}
\newcommand{\C}{\mathbb{C}}
\newcommand{\T}{\mathbb{T}}

\newcommand{\E}{\mathbb{E}}

\newcommand{\lcm}{\mathrm{lcm}}

\newcommand{\eps}{\varepsilon}

\newcommand{\CB}{\mathcal{B}}
\newcommand{\CI}{\mathcal{I}}
\newcommand{\str}{\mathrm{str}}

\newcommand{\unf}{\mathrm{unf}}

\makeatletter
\let\@@pmod\pmod
\DeclareRobustCommand{\pmod}{\@ifstar\@pmods\@@pmod}
\def\@pmods#1{\mkern4mu({\operator@font mod}\mkern 6mu#1)}
\makeatother

\renewcommand{\leq}{\leqslant}
\renewcommand{\geq}{\geqslant}

\newtheorem{theorem}{Theorem}[section]

\newtheorem{proposition}[theorem]{Proposition}
\newtheorem{lemma}[theorem]{Lemma}

\theoremstyle{definition}
\newtheorem{definition}[theorem]{Definition}
\newtheorem*{remark}{Remark}

\usepackage[normalem]{ulem}

\usepackage[makeroom]{cancel}

\numberwithin{theorem}{section}

\title{Quantitative bounds in a popular polynomial Szemer\'{e}di theorem}

\author{Xuancheng Shao}
\address{Department of Mathematics, University of Kentucky, Lexington, KY, 40506, USA}
\email{xuancheng.shao@uky.edu}

\author{Mengdi Wang}
\address{\'{E}cole Polytechnique F\'{e}d\'{e}rale de Lausanne (EPFL), Lausanne, Switzerland}
\email{mengdi.wang@epfl.ch}

\begin{document}

\begin{abstract}
We obtain polylogarithmic bounds in the polynomial Szemer\'{e}di theorem when the polynomials have distinct degrees and zero constant terms. Specifically, let $P_1, \dots, P_m \in \mathbb Z[y]$ be polynomials with distinct degrees, each having zero constant term. Then there exists a constant $c = c(P_1,\dots,P_m) > 0$ such that any subset $A \subset \{1,2,\dots,N\}$ of density at least $(\log N)^{-c}$ contains a nontrivial polynomial progression of the form $x, x+P_1(y), \dots, x+P_m(y)$. In addition, we prove an effective ``popular'' version, showing that every dense subset $A$ has some non-zero $y$ such that the number of polynomial progressions in $A$ with this difference $y$ is asymptotically at least as large as in a random set of the same density as $A$.
\end{abstract}

\maketitle

\section{Introduction}

The quest to understand the ubiquity of arithmetic structures within subsets of integers has driven some of the most profound advances in additive combinatorics. Szemer\'{e}di's landmark theorem states that every subset of the integers with positive density contains non-trivial $k$-term arithmetic progressions ($k$-APs) for every natural number $k$. A pivotal contribution came from Gowers \cite{Go}, who introduced uniformity norms to furnish a quantitative bound for Szemer\'{e}di's theorem, establishing that any subset of $\{1,2,\dots,N\}$ which does not contain non-trivial $k$-APs must have density $O((\log\log N)^{-c_k})$ for some constant $c_k>0$. Over the years, the pursuit of quantitative bounds in Szemer\'{e}di's theorem has remained extremely active,  leading to many significant developments. For example, Kelley and Meka \cite{kelley-meka} established a bound of the shape $O(\exp(-c(\log N)^{1/12}))$ for the case $k=3$, where the exponent $1/12$ is later improved to $1/9$ by Bloom and Sisask \cite{bloom-sisask}. Green and Tao \cite{GT-poly-4aps} achieved a polylogarithmic bound for $k=4$. Leng, Sah, and Sawhney \cite{LSS2} improved Gowers' bound to $O(\exp(-(\log \log N)^{c_k}))$ for general $k\geq 5$.

The scope of such results broadened significantly with the polynomial Szemer\'{e}di theorem of Bergelson and Leibman \cite{BL}, which established that every subset of the integers with positive density contains polynomial progressions of the form $x, x+P_1(y), \dots, x+P_m(y)$ with $y \neq 0$, where $P_1,\dots,P_m \in \Z[y]$ and each has zero constant term. Their proof was rooted in ergodic theory and was thus non-quantitative.

The \emph{density increment method} has long served as the primary tool for deriving quantitative bounds for general additive patterns. However, applying this strategy to more complex polynomial configurations introduces additional difficulties; see \cite[Section 2.2]{Pre14} for a detailed discussion.   While a quantitative version of the general polynomial Szemer\'edi theorem  remains an open challenge, progress has been made in certain special cases.   We summarize some of these developments below and refer the reader to \cite{Pe} for  for a more in-depth overview of the history of quantitative results in the polynomial Szemer\'{e}di theorem both in the integers and in the finite field setting.

\begin{itemize}
\item For Sark\"{o}zy's theorem on square differences (i.e. $x, x+y^2$), Green and  Sawhney \cite{GS-fs} achieved quantitative bounds of the form $O(\exp(-c\sqrt{\log N}))$.

\item When $P_1,\cdots,P_m$ are homogeneous and have the same degree, Prendiville \cite{Pre14} achieved bounds of the form $O((\log\log N)^{-c})$ (which can likely be improved to $O(\exp(-(\log \log N)^{c}))$ using the quasipolynomial bounds in the inverse theorem for the Gowers norms from \cite{LSS}, together with the density increment strategy of Heath-Brown and Szemer\'edi).

\item For the nonlinear Roth configuration (i.e. $x, x+y, x+y^2$), Peluse and Prendiville \cite{PP1,PP2} achieved quantitative bounds of the form $O((\log N)^{-c})$.

\item Peluse \cite{Pe} extended these techniques from the nonlinear Roth configuration to more general polynomial progressions when $P_1,\dots, P_m$ have distinct degrees and each has zero constant term. She obtained bounds of the form $O((\log\log N)^{-c})$ for subsets of $\{1,2,\dots,N\}$ that do not contain such configurations.
\end{itemize}

The first main result of this paper is an improvement of the quantitative bound in Peluse's result \cite[Theorem 1.1]{Pe}.

\begin{theorem}[Density bound]\label{thm:density}
Let $P_1,\dots,P_m \in \Z[y]$ be polynomials with distinct degrees, each having zero constant term. If $A \subset \{1,2,\dots,N\}$ does not contain a polynomial progression of the form
$$
x, x+P_1(y), \dots, x+P_m(y) \ \ \text{with }y \neq 0,
$$
then 
$$
|A| \ll \frac{N}{(\log N)^c}
$$
for some constant $c>0$ depending only on $P_1,\dots,P_m$.
\end{theorem}

In the ``structure vs. randomness" framework underlying the density increment argument, if a set $A$ is not pseudorandom, then the number of additive configurations in this set $A$ may deviate significantly from the random bound. Bergelson, Host, and Kra \cite{BHK} posed the question of whether there always exists a nonzero difference  $y$ for which the number of configurations in $A$ with  difference parameter $y$ is at least as large as the random bound. For example, in the case when the configurations are $3$-APs, the question asks for $A \subset \{1,2,\cdots,N\}$ with $|A| = \delta N$, whether there exists $y \neq 0$ such that
$$
\#\{x \in A: x, x+y, x+2y \in A\} \geq (\delta^3 - o(1))N.
$$
Such $y$ is loosely referred as a ``popular difference". This topic in the ergodic setting was initially explored by Bergelson, Host, and Kra \cite{BHK} and also studied in \cite{Fra, FK, DLMS}. We list some notable results below in the arithmetic setting, and refer the readers to \cite{PPS} for more contexts surrounding the popular difference results.

\begin{itemize}
\item The existence of popular differences is proved by Green \cite{G-regularity} for $3$-APs and by Green and Tao \cite{GT-regularity} for $4$-APs. When $k\geq 5$, a counterexample  was built by Ruzsa in the appendix to \cite{BHK}.
\item The quantitative aspects of popular common differences in 3-APs were studied in \cite{FP, FP2, FPZ}.\item Linear configurations were classified in \cite{SSZ} according to when popular differences are guaranteed to exist.
\item Popular differences for corners $(x_1,x_2), (x_1+y,x_2), (x_1,x_2+y)$ are studied in \cite{Mandache, FSSSZ}.
\item The existence of popular differences for special polynomial progressions is proved by Lyall and Magyar  \cite{LM} for the square differences $x, x+y^2$ and by Peluse, Prendiville, and the first author \cite{PPS} for the nonlinear Roth pattern $x, x+y, x+y^2$.
\end{itemize}

In this paper, we prove, with effective bounds, the existence of popular differences in general polynomial progressions when the polynomials have distinct degrees and each has zero constant term.

\begin{theorem}[Popular difference]\label{thm:popular}
Let $P_1,\dots,P_m \in \Z[y]$ be polynomials with distinct degrees, each having zero constant term. Let $A \subset \{1,2,\dots,N\}$ be a subset with $|A| \geq \delta N$ and let $\eps \in (0,1/2)$. Then either $N \leq \exp(\exp(\eps^{-O_{P_1,\dots,P_m}(1)}))$ or there exists a positive integer $y$ such that
$$
\#\{x \in A: x+P_i(y) \in A\text{ for each }1 \leq i \leq m\} \geq (\delta^{m+1}-\eps)N.
$$
\end{theorem}

The proofs of Theorems \ref{thm:density} and \ref{thm:popular} are based on an inverse theorem. As the final part of the introduction, we present a simplified version of our more general result, Theorem \ref{thm:inverse}. Suppose that $N, M \geq 1$,  $P_1,\dots,P_m \in \Z[y]$ are polynomials and $f_0,\dots,f_m:\Z\rightarrow\C$ are functions  supported on the interval $[N]$. We define the associated counting operator as follows:
\begin{align}\label{counting-mn}
\Lambda_{P_1,\dots,P_m}^{N,M} (f_0,\dots, f_m) := \frac{1}{NM} \sum_{x \in \Z} \sum_{y \in [M]} f_0(x) f_1(x+P_1(y)) \cdots f_m(x+P_m(y)).
\end{align}

%{\color{violet}
%\begin{theorem}[Inverse theorem]
%Let $P_1,\dots,P_m \in \Z[y]$ be polynomials with $(C,q)$-coefficients such that they have distinct degrees. Let $d = max_{1 \leq i \leq m}\deg P_i$. Let $N, M$ be positive integers with $M \leq (N/q^{d-1})^{1/d}$, and let $\delta \in (0,1/2)$. Let $f_0,\dots,f_m:\Z\rightarrow \C$ be $1$-bounded functions supported on $[N]$. Suppose that
%$$
%\bigabs{\Lambda_{P_1,\dots,P_m}^{N,M}(f_0,\dots,f_m)} \geq \delta.
%$$
%Then either $M \ll_{C,d} (q/\delta)^{O_d(1)}$ or there exist positive integers $q'\ll_{C,d}\delta^{-O_d(1)}$ and $b = O_d(1)$, as well as a $1$-bounded function $\phi_i:\Z\rightarrow\C$ for each $1 \leq i \leq m$ which is $O_{C,d}((q/\delta)^{O_d(1)}M^{-\deg P_i})$-Lipschitz along $q'q^b\cdot\Z$, such that
%$$
%\bigabs{\sum_{x \in \Z}f_i(x)\phi_i(x)}\gg_{C,d} \delta^{O_d(1)}N.
%$$
%\end{theorem}}

\begin{theorem}[Inverse theorem]\label{thm:inverse-simple}
Let $P_1,\dots,P_m \in \Z[y]$ be polynomials with distinct degrees, and let their coefficients be bounded by $C$. Let $d = max_{1 \leq i \leq m}\deg P_i$. Let $N, M$ be positive integers with $M \asymp_C N^{1/d}$, and let $\delta \in (0,1/2)$. Suppose that $f_0,\dots,f_m:\Z\rightarrow \C$ are $1$-bounded functions supported on $[N]$. If
$$
\bigabs{\Lambda_{P_1,\dots,P_m}^{N,M}(f_0,\dots,f_m)} \geq \delta.
$$
Then either $M \ll_{C,d} \delta^{-O_d(1)}$ or there exists a positive integer $q\ll_{C,d}\delta^{-O_d(1)}$ and a $1$-bounded function $\phi_i:\Z\rightarrow\C$ for each $1 \leq i \leq m$ which is $O_{C,d}(\delta^{-O_d(1)}M^{-\deg P_i})$-Lipschitz along $q\cdot\Z$, such that
$$
\bigabs{\sum_{x \in \Z}f_i(x)\phi_i(x)}\gg_{C,d} \delta^{O_d(1)}N.
$$
\end{theorem}

The definition of $C$-Lipschitz function is given in Definition \ref{def-lip}. For intuition, we encourage the reader to think of $\phi_i$ as a phase function of the form $x \mapsto e(\alpha_ix)$ with $\norm{q\alpha_i}\ll_{C,d}\delta^{-O_d(1)}M^{-\deg P_i}$. Alternatively, $\phi_i$ is almost constant on progressions of length $\Theta(\delta^{O_d(1)}M^{\deg P_i})$ and common difference $q$. Even though $q$ is independent of $i$, this is unimportant since one could simply take $q$ to be the least common multiple of all $q_i$'s if $q$ is allowed to depend on $i$, while maintaining the upper bound for $q$.

Compared to Peluse's inverse theorem \cite[Theorem 3.3]{Pe} which concludes correlation  only for the $f_i$ corresponding to the lowest degree polynomial, our inverse theorem gives such conclusions for every $f_i$. This generalization, which is fundamental in getting the polylogarithmic bound in Theorem \ref{thm:density} and also in obtaining the popular difference result, is largely based on arguments of Peluse \cite{Pe} incorporated with some technical generalizations. Instead of proving Theorem \ref{thm:inverse-simple} directly, we will in fact prove a slightly more general version of it, Theorem \ref{thm:inverse}.
See the discussions following the statement of Theorem \ref{thm:inverse} for a brief overview of  the proof ideas.

We would like to make some remarks on general polynomial sequences. When the leading coefficients of the highest-degree terms among $P_1,\dots, P_m$ are distinct, the sequence can be controlled by a Gowers $U^s$-norm for some $s>2$. In other cases, it can only be controlled by an average of box norms. Ideally, in such situations, some suitable local inverse theorems can be applied. As suggested by Leng--Sah--Sawhney \cite{LSS2}, one can obtain quantitative density bounds for linear forms, since the inverse theorem implies that $f$ correlates with a periodic function $\phi$. However, the period of $\phi$ is typically of order $N^{1-\eta}$ for some small $\eta>0$, which means that the density increment achieved at each step is not sufficiently effective. This is also the reason why the density increment argument fails for general polynomial sequences.

The rest of the paper is organized as follows. In Section \ref{sec:inverse} we prove our main inverse theorem, Theorem \ref{thm:inverse} (of which Theorem \ref{thm:inverse-simple} is a special case). In Section \ref{sec:factors} we develop the theory of \emph{local factors} which are necessary for our weak regularity lemmas. In Section \ref{sec:density-increment} we prove Theorem \ref{thm:density} by establishing a regularity lemma and using the density increment argument, generalizing Peluse and Prendiville's work \cite{PP2} on the nonlinear Roth pattern. In Section \ref{sec:popular} we prove Theorem \ref{thm:popular} by establishing a second weak regularity lemma.

\subsection*{Acknowledgements}

The authors  thank Sean Prendiville for valuable discussions, particularly regarding an earlier version of the proof of Theorem \ref{thm:popular}, and for his encouragement. 
We are grateful to the anonymous referee for
a careful reading of the paper and for numerous helpful comments and corrections.
XS was supported by NSF grant DMS-2200565.
MW was partially supported by the Swiss National Science Foundation grant TMSGI2-2112.

%%%%%%%%%%%%%%%%%%%%%%%%%%%%%%%%%%%%%%%%%%%%%%%%%%%%%%%%%%%%%%%%%%%%%%%%%%%%%%%%%%%%%%%%%%%%%%%%%%%%%%%%%%%%%%%%%%%%%%%%%%%%%%

\section{Inverse theorems}\label{sec:inverse}

The goal of this section is to prove the inverse theorem which is the key to both Theorems \ref{thm:density} and \ref{thm:popular}. It asserts that if the counting expression (\ref{counting-mn}) is large, then the  functions exhibit correlation with certain Lipschitz functions. We adopt the same definition as \cite[Definition 6.1]{PP2}.

\begin{definition}[$C$-Lipschitz] \label{def-lip}
A function $\phi:\Z\to\C$ is said to be $C$--\emph{Lipschitz} along $q\cdot\Z$ if for all $x,y\in \Z$ we have
\[
|\phi(x+qy)-\phi(x)|\leq C|y|.
\]
	
\end{definition}

Note that if  $\alpha\in\T$ is a frequency such that $\norm{q\alpha}$ is sufficiently small for some positive integer $q\geq 1$, then the linear phase function $e(\cdot\alpha)$ is $(2\pi \norm{q\alpha})$-Lipschitz along $q\cdot\Z$. Thus Lipschitz functions serve as an analogue to  linear phase functions. The following lemma  provides another class of Lipschitz functions.

\begin{lemma}\label{lipschitz}
	Let $q,H$ be positive integers and let $f:\Z\to\C$ be a 1-bounded function. Then the function
	\[
	\phi(x) = \E_{h_1,h_2\in[H]} f(x+q(h_1-h_2))
	\]
	is $O(H^{-1})$-Lipschitz along $q\cdot\Z$.
\end{lemma}
\begin{proof}
\cite[Lemma 6.2]{PP2}.	
\end{proof}

 Before stating our general inverse theorems, we recall the notation of $(C,q)$-coefficients from \cite[Definition 3.1]{Pe}, which serves as an important index for tracking iteration steps.

\begin{definition}\label{coefficient}
We say a polynomial $P(y)=a_dy^d+\cdots+a_1y \in \Z[y]$ has \emph{$(C,q)$-coefficients} if $|a_i|\leq C|a_d|$ for all $1\leq i \leq d-1$ and $a_d=a_d'q^{d-1}$ for some $a_d' \in \Z$ with $0< |a_d'|\leq C$.	
\end{definition}

We are now ready to state the inverse theorem, of which Theorem \ref{thm:inverse-simple} from the introduction is a special case corresponding to the situation where the polynomials have $(C,1)$-coefficients and $M$ is relatively large. Recall the counting operator defined in \eqref{counting-mn}.

\begin{theorem}[Inverse theorem]\label{thm:inverse}
Let $P_1,\dots,P_m \in \Z[y]$ be polynomials with $(C,q)$-coefficients, and assume they have distinct degrees. Let $d = max_{1 \leq i \leq m}\deg P_i$. Let $N, M$ be positive integers with $M \leq (N/q^{d-1})^{1/d}$, and let $\delta \in (0,1/2)$. Suppose that $f_0,\dots,f_m:\Z\rightarrow \C$ are $1$-bounded functions supported on $[N]$. Let
$$
\bigabs{\Lambda_{P_1,\dots,P_m}^{N,M}(f_0,\dots,f_m)} \geq \delta.
$$
Then either $M \ll_{C,d} (q/\delta)^{O_d(1)}$, or there exist positive integers $q'\ll_{C,d}\delta^{-O_d(1)}$ and $b = O_d(1)$, as well as a $1$-bounded function $\phi_i:\Z\rightarrow\C$ for each $1 \leq i \leq m$ which is $O_{C,d}((q/\delta)^{O_d(1)}M^{-\deg P_i})$-Lipschitz along $q'q^b\cdot\Z$, such that
$$
\bigabs{\sum_{x \in \Z}f_i(x)\phi_i(x)}\gg_{C,d} \delta^{O_d(1)}N.
$$
\end{theorem}

Assuming $\deg P_1=\min_{1\leq i\leq m}\deg P_i$ without loss of generality, we note that Peluse's inverse theorem \cite[Theorem 3.3]{Pe} implies the correlation result only for the function $f_1$ in Theorem \ref{thm:inverse}, under the assumption that $M = (N/q^{d-1})^{1/d}$ (or, slightly more generally, $M \asymp_C (N/q^{d-1})^{1/d}$). The new ingredient in our proof of Theorem \ref{thm:inverse} involves a crucial inductive step which passes from the correlation conditions for $f_1,\cdots,f_{j-1}$ to the correlation condition for $f_j$. Loosely speaking, if $f_1,\dots,f_{j-1}$ correlates with Lipschitz functions $\phi_1,\dots,\phi_{j-1}$, respectively, then we will deduce that
$$
\bigabs{ \Lambda_{P_1,\dots,P_m}^{N,M}(f_0,\phi_1,\dots,\phi_{j-1},f_j,\dots,f_m) } \gg \delta^{O(1)}.
$$
Assume, for simplicity, that $\phi_1,\dots,\phi_{j-1}$ are constant functions, this inequality then implies that there is a 1-bounded function $g_0$ such that
\[
\bigabs{ \Lambda_{P_j,\dots,P_m}^{N,M}(g_0,f_j,\dots,f_m) } \gg \delta^{O(1)},
\]
from which we will then deduce that $f_j$ also correlates with a Lipschitz function. This idea was discussed in \cite[Section 7]{PP2}, and our implementation of it involves technical generalizations of Peluse's arguments \cite{Pe}.

These arguments outlined above will be carried out in Sections \ref{sec:inverse1} and \ref{sec:inverse2}, proving Theorem \ref{thm:inverse} when $M \asymp_C (N/q^{d-1})^{1/d}$. In Section \ref{sec:inverse3}, we will generalize to all $M \leq (N/q^{d-1})^{1/d}$ by working in subintervals of $[N]$ of appropriate lengths.

\subsection{Inductive step}\label{sec:inverse1}

In this subsection, we prove Proposition \ref{partial-ii} as a necessary step allowing us to pass from the correlation conditions for $f_1,\dots,f_{j-1}$ to the correlation conclusion for $f_j$.

\begin{proposition}[Partial correlation]\label{partial-ii}
Let $N,M>0$ be sufficiently large numbers and $q\in\N$. Let $P_1,\dots,P_m\in\Z[y]$ be polynomials with $(C,q)$-coefficients such that $\deg P_1< \cdots  <\deg P_m$. Define $d=\deg P_m$ and assume $1/C \leq q^{d-1}M^d/N \leq C$. Let $\delta \in (0,1/2)$. We have either $N\ll_{C,d} (q/\delta)^{O_{d}(1)}$ or the following conclusions hold.

%	\item  If $f_0,\dots,f_m:\Z\to\C$ are 1-bounded functions supported on the interval $[N]$ and
%	\[\bigabs{ \Lambda_{P_1,\dots,P_m}^{N,M}(f_0,\dots,f_m) } \geq  \delta,
%	\]
%	 then  there exist positive integers $q' \ll_{C,d} \delta^{-O_{d}(1)}$, $b \ll_{d} 1$ along with a 1-bounded function $\phi_1:\Z\to\C$ that is $O_{C,d}((q/\delta)^{O_{d}(1)} M^{-\deg P_1})$-Lipschitz along $q'q^b\cdot\Z$, such that
%\[
%\bigabs{ \sum_x f_1(x) \phi_1(x) }\gg_{C,d} \delta^{O_{d}(1)} N.
%\]	

Suppose that $1\leq j\leq m$, and for each $1\leq i<j$,  let $\phi_i:\Z\to\C$ be a 1-bounded $O_{C,d}((q/\delta)^{O_{d}(1)} M^{-\deg P_i})$-Lipschitz function along $q_iq^{b_i}\cdot\Z$ for some integers $q_i\ll_{C,d}\delta^{-O_{d}(1)}$ and $b_i\ll_{d} 1$. Let $f_0,f_j,\cdots,f_m:\Z\to\C$ be arbitrary 1-bounded functions supported on $[N]$. If
 \begin{align}\label{ineq-j}
\bigabs{ \Lambda_{P_1,\dots,P_m}^{N,M}(f_0,\phi_1,\dots,\phi_{j-1},f_j,\dots,f_m) } \geq  \delta,
\end{align}
then there exist positive integers $q_j \ll_{C,d} \delta^{-O_{d}(1)}$, $b_j \ll_{d} 1$ and a 1-bounded function $\phi_j:\Z\to\C$ which is $O_{C,d}((q/\delta)^{O_{d}(1)} M^{-\deg P_j})$-Lipschitz along $q_jq^{b_j}\cdot\Z$, such that
\[
\bigabs{ \sum_x f_j(x) \phi_j(x) }\gg_{C,d} \delta^{O_{d}(1)} N.
\]	
\end{proposition}

This is an extension of \cite[Theorem 3.3]{Pe} which essentially corresponds to the $j=1$ case. The assumption $q^{d-1}M^d/N = 1$ is assumed in \cite[Theorem 3.3]{Pe}, but it can easily be relaxed to $q^{d-1}M^d/N \asymp 1$ by following the proof.

The proof of Proposition \ref{partial-ii} adapts the proof of \cite[Theorem 3.3]{Pe} in \cite[Section 9]{Pe}. We outline the arguments before presenting the details. For $N, M \geq 1$, polynomials $P_1,\cdots,P_m \in \Z[y]$, functions $f_0,\cdots,f_{\ell}:\Z\rightarrow\C$ supported on $[N]$ and characters $\psi_{\ell+1},\cdots,\psi_m:\Z\rightarrow\C$ of the form $\psi_i(x) = e(\alpha_ix)$ with $\alpha_i\in\R$, we define the counting operator
\begin{equation*}
\begin{split}
& \Lambda_{P_1,\cdots,P_m}^{N,M}(f_0,\cdots,f_{\ell};\psi_{\ell+1},\cdots,\psi_m) \\
&  := \frac{1}{NM} \sum_{x \in \Z} \sum_{y \in [M]} f_0(x) f_1(x+P_1(y)) \cdots f_{\ell}(x+P_{\ell}(y)) \psi_{\ell+1}(P_{\ell+1}(y)) \cdots \psi_m(P_m(y)).
\end{split}
\end{equation*}
When $\ell=m$ this notation becomes \eqref{counting-mn}. The first step in the proof of Proposition \ref{partial-ii} involves repeated applications of the following lemma which is \cite[Lemma 3.11]{Pe}.

\begin{lemma}\label{peluse-lemma3.11}
Let $N,M>0$ be sufficiently large numbers, $q\in\N$, and $2 \leq \ell \leq m$. Let $P_1,\dots,P_m\in\Z[y]$ be polynomials such that $P_1,\cdots,P_{\ell}$ have $(C,q)$-coefficients and $\deg P_1< \cdots  <\deg P_m$. Let $d = \deg P_{\ell}$ and let $c_{\ell}$ be the leading coefficient of $P_{\ell}$. Assume $1/C \leq q^{d-1}M^d/N \leq C$. Let $\delta \in (0,1/2)$. We have either $N\ll_{C,d} (q/\delta)^{O_{d}(1)}$ or the following conclusions hold.

Let $f_0,\cdots,f_{\ell}: \Z\rightarrow \C$ be $1$-bounded functions supported on $[N]$ and let $\psi_{\ell+1},\cdots,\psi_m$ be characters. If
$$
\bigabs{\Lambda_{P_1,\cdots,P_m}^{N,M}(f_0,\cdots,f_{\ell};\psi_{\ell+1},\cdots,\psi_m)} \geq \delta,
$$
then we have
$$
\E_{\substack{u, h = 0,\cdots,|c'|-1 \\ 0 \leq w < N/(c'C'N')}} \bigabs{\Lambda_{P_1^h,\cdots,P_m^h}^{C'N',M'}(f_0^{u,h,w},\cdots,f_{\ell-1}^{u,h,w}; \psi_{\ell,u}, \psi_{\ell+1},\cdots,\psi_m)} \gg_{C,d} \delta^{O_d(1)} 
$$
for some characters $\psi_{\ell,u}$, where $C' \asymp_d C$, $c' := d!c_{\ell}$, $M' := M/|c'|$, $N' := M'^{\deg P_{\ell-1}} (q|c'|)^{\deg P_{\ell-1}-1}$,
$$
P_i^h(z) := \begin{cases} \frac{1}{c'}(P_i(c'z+h) - P_i(h)) & 1 \leq i \leq \ell-1, \\ P_i(c'z+h) - P_i(h) & \ell \leq i \leq m, \end{cases}
$$
and
$$
f_i^{u,h,w}(x) := \begin{cases} (f_0\psi_{\ell,u})\Big(c'x + c'C'N'w - P_{\ell}(h) - u\Big) \cdot 1_{[C'N']}(x) & i=0, \\
f_i\Big(c'x + c'C'N'w + P_i(h) - P_{\ell}(h)-u\Big) \cdot 1_{[C'N']}(x) & 1 \leq i \leq \ell-1. \end{cases}
$$
\end{lemma}

This allows us to essentially replace $f_{j+1},\cdots,f_m$ in \eqref{ineq-j} by characters. By applying Weyl's exponential sum estimates stated below (see \cite[Lemma 1.1.16]{Tao}), one can deduce that these characters must be in major arcs.

\begin{lemma}\label{vinogradov-lemma}
Let $N\geq 1$ and $0<\delta<1$. Suppose that $\alpha_1,\dots,\alpha_d\in\R $ and 
\[
\Bigabs{ \sum_{n\leq N} e\bigbrac{ \alpha_d n^d +\cdots + \alpha_1n}} \geq\delta N,
\]	
then there exists a positive integer $1\leq q\ll \delta^{-O_d(1)}$ such that $\norm{q\alpha_i}\ll\delta^{-O_d(1)}N^{-i}$ for each $1\leq i\leq d$.
\end{lemma}

At this stage, we essentially have
$$
\bigabs{ \Lambda_{P_1,\dots,P_m}^{N,M}(f_0,\phi_1,\dots,\phi_{j-1},f_j,\psi_{j+1},\cdots,\psi_m) } \gg \delta^{O(1)},
$$
where $\psi_{j+1},\cdots,\psi_m$ are major-arc linear phase functions. By the Lipschitz properties of $\phi_1,\cdots,\phi_{j-1}$, we may split the range of $y \in [M]$ into arithmetic progressions of length $|P| \approx (\delta/q)^{O(1)}M$ such that if $y \in P$ then
$$
\phi_i(x) \approx \phi_i(x + P_i(y)) \text{ for }1 \leq i \leq j-1, \ \ \psi_i(x) \approx \psi_i(x+P_i(y)) \text{ for } j+1 \leq i \leq m.
$$
This leads to
$$
\frac{1}{N}\sum_{x \in \Z} \sum_{y \in P} f_0(x) f_j(x + P_j(y)) \gg \delta^{O(1)}|P|
$$
(after modifying $f_0$ appropriately). The desired Lipschitz property for $f_j$ follows by applying Cauchy-Schwarz and van der Corput's inequality as in \cite[Lemma 4.2]{Pe}.

\begin{proof}[Proof of Proposition \ref{partial-ii}]
Assume that the inequality (\ref{ineq-j}) holds. Following the proof of \cite[Theorem 3.3]{Pe}, applying Lemma \ref{peluse-lemma3.11} a total of $m-j$ times yields
\begin{multline} \label{apply-pe3.11}
	\E_{\substack{u_i,h_i=0,\dots,|c_i|-1 \\
0\leq w_i<(C_{i+1}N_{i+1}/|c_i|)/C_iN_i \\
j+1\leq i\leq m}}   \\  \Bigabs{ \Lambda_{ P_1^{\underline{h}},\dots, P_m^{\underline{h}} }^{C_{j+1}N_{j+1}, M_{j+1} } (f_0^{\underline u, \underline h, \underline w}, \phi_1^{\underline u, \underline h, \underline w} ,\dots, \phi_{j-1}^{\underline u, \underline h, \underline w}, f_j^{\underline u, \underline h, \underline w};\psi_{j+1}^{\underline u, \underline h, \underline w},\dots,\psi_m^{\underline u, \underline h, \underline w})} \gg_{C,d}  \delta^{O_{d} (1)}.
\end{multline}
The parameters involved in the equality satisfy the following conditions: $C_{m+1}=1$, $N_{m+1}=N$, $b_i\ll_{d}1$, $c_i\asymp_{C, d} q^{b_i}$,  $M_i=M/\prod_{i\leq j\leq m}|c_j|$, $C_i\asymp_{C,d} 1$ and $N_i = M_i^{\deg P_i-1}(q|c_i\cdots c_m|)^{\deg P_{i-1}-1}$ for each $i=j+1,\dots,m$. Besides, $f_0^{\underline u, \underline h, \underline w}$ is a 1-bounded function,
\begin{align}\label{phi-i}
 \phi_i^{\underline u, \underline h, \underline w} (x) =\phi_i(c_{j+1}\cdots c_m x +Q_i(\underline u, \underline h, \underline w))
\end{align}
for $1 \leq i \leq j-1$, where $Q_i(\underline u, \underline h, \underline w)$ is a function in terms of $\underline u, \underline h, \underline w$ and $P_i^{\underline h}$, and $f_j^{\underline u, \underline h, \underline w}(x)$ equals $1_{[C_{j+1}N_{j+1}]}(x)$ times
\[
f_j \Bigbrac{c_{j+1}\cdots c_m x + \sum_{j< i\leq m} (c_{i+1}\cdots c_m)[w_ic_iC_i N_i - u_i + [P_j^{h_m,\dots,h_{i+1}} (h_i) - P_i^{h_m,\dots,h_{i+1}} (h_i)]]}.
\]

Next we will remove the functions $\phi_i^{\underline u, \underline h, \underline w}$ and $\psi_i^{\underline u, \underline h, \underline w}$ from (\ref{apply-pe3.11}) step by step. First,  for each character $\psi_i^{\underline u, \underline h, \underline w}$ ($j+1 \leq i \leq m$), we can express it as $\psi_i^{\underline u, \underline h, \underline w} (x) =e (\beta_i^{\underline u, \underline h, \underline w}x)$ for some frequency $\beta_i^{\underline u, \underline h, \underline w}\in\R $. By applying the Cauchy-Schwarz inequality $O_{\deg P_1,\dots,\deg P_j} (1)$ times and then invoking Lemma \ref{vinogradov-lemma}, we may find positive integers  $t'\ll_{C,d}\delta^{-O_{d}(1)}$ and $b_i\ll_{d}1$ such that
\[
\norm{ t' q^{b_i}\beta_i^{\underline u, \underline h, \underline w} } \ll_{C,d} \delta^{-O_{d}(1)} /M_{j+1}^{\deg P_i}
\]
for each $j+1\leq i\leq m$. The range of summation $y \in [M_{j+1}]$ in \eqref{apply-pe3.11} can be split into arithmetic progressions, each of which has length $M_{j+1}^* =\Theta_{C,d} ((\delta/q)^{O_{d}(1)} M_{j+1})$ and  common difference $t'q^s$ for some $s\ll_{d}1$,  such that  $\psi_i^{\underline u, \underline h, \underline w}$ is nearly constant within each arithmetic progression. Moreover, the number of such arithmetic progressions is $O(t'q^sM_{j+1}/M_{j+1}^*)$. Thus, after a change of variables, we deduce from inequality (\ref{apply-pe3.11})  that
\begin{multline}\label{eq2.4}
\delta^{O_{d} (1)} \ll_{C,d} \E_{ \substack{u_i,h_i=0,\dots,|c_i|-1\\
0\leq w_i<(C_{i+1}N_{i+1}/|c_i|)/C_iN_i \\
j+1\leq i\leq m\\
k_{\underline u, \underline h, \underline w} \in [M_{j+1}/M_{j+1}^*]\\
k'_{\underline u, \underline h, \underline w} \in [t'q^s]}}	
\Bigabs{(C_{j+1}N_{j+1})^{-1}\sum_x \E_{y\in[M_{j+1}^*]} f_0^{\underline u, \underline h, \underline w} (x) \\
\phi_1^{\underline u, \underline h, \underline w} (x+P_1^{\underline h}(t'q^s (y-M_{j+1}^* k'_{\underline u, \underline h, \underline w}) - k_{\underline u, \underline h, \underline w})) \cdots f_j^{\underline u, \underline h, \underline w} (x + P_j^{\underline h} (t'q^s (y-M_{j+1}^* k'_{\underline u, \underline h, \underline w}) - k_{\underline u, \underline h, \underline w})) }.
\end{multline}

Next, from the Lipschitz property of $\phi_i$ and \eqref{phi-i}, one can deduce that, for each $1 \leq i \leq j-1$, the function $\phi_i^{\underline u, \underline h, \underline w}$  is  $O_{C,d}((q/\delta)^{O_{d}(1)} M^{-\deg P_i})$-Lipschitz  along $q_iq^{b_i}\cdot\Z$, noting that $|c_{j+1}\cdots c_m|\ll_{C,d} q^{O_d(1)}$. 
By increasing $t'$ and $s$ if necessary, we may ensure that $\lcm[q_1,\cdots,q_{j-1}] \mid t'$ and $s \geq \max_{1 \leq i \leq j-1}b_i$, while maintaining the bounds $t' \ll_{C,d} \delta^{-O_d(1)}$ and $s\ll_d1$. Thus $q_iq^{b_i} \mid t'q^s$ for each $1 \leq i \leq j-1$. It follows that the dependence on $y$ in each of the terms involving $\phi_1^{\underline u, \underline h, \underline w}, \cdots, \phi_{j-1}^{\underline u, \underline h, \underline w}$ can be dropped at a negligible error. More specifically, writing $Q(y): = P_i^{\underline h}(t'q^s(y-M_{j+1}^*k_{\underline u, \underline h, \underline w}')-k_{\underline u, \underline h, \underline w})$, we have for $y \in [M_{j+1}^*]$,
$$
\bigabs{ \phi_i^{\underline u, \underline h, \underline w} (x+Q(y)) - \phi_i^{\underline u, \underline h, \underline w} (x+Q(0))}
\ll_{C,d} (q/\delta)^{O_{d} (1)} M^{-\deg P_i} |Q(y) - Q(0)| 
$$
since $q_iq^{b_i} \mid Q(y)-Q(0)$. Since $P_i^{\underline h}$ has degree $\deg P_i$ and all of its coefficients are $\ll_{C,d} q^{O_d(1)}$, it follows that 
$$
|Q(y) - Q(0)| \ll_{C,d} (q/\delta)^{O_d(1)} M^{\deg P_i-1} \cdot M_{j+1}^*.
$$
Hence, by decreasing $M_{j+1}^*$ if necessary (while still maintaining the bound $M_{j+1}^* \gg_{C,d} (\delta/q)^{O_d(1)}M_{j+1}$), we may ensure that the function $\phi_i^{\underline u, \underline h, \underline w} (x+Q(y))$ can be replaced by $\phi_i^{\underline u, \underline h, \underline w} (x+Q(0))$ in \eqref{eq2.4}.
Therefore, we conclude that
\begin{multline*}
\delta^{O_{d} (1)} \ll_{C,d} \E_{ \substack{u_i,h_i=0,\dots,|c_i|-1\\
0\leq w_i<(C_{i+1}N_{i+1}/|c_i|)/C_iN_i \\
j+1\leq i\leq m\\
k_{\underline u, \underline h, \underline w} \in [M_{j+1}/M^*_{j+1}]\\
k'_{\underline u, \underline h, \underline w} \in [t'q^s]}}  (C_{j+1}N_{j+1})^{-1}\\
\Bigabs{\sum_x \E_{y\in[M^*_{j+1}]} g_0^{\underline u, \underline h, \underline w} (x) 
f_j^{\underline u, \underline h, \underline w} (x + P_j^{\underline h} (t'q^s (y-M^*_{j+1} k'_{\underline u, \underline h, \underline w}) - k_{\underline u, \underline h, \underline w})) },
\end{multline*} 
 where $ g_0^{\underline u, \underline h, \underline w} (x)$ is $1$-bounded.
 This inequality has the same shape as the fourth from last inequality on  \cite[Page 52]{Pe} with $f_0,f_1$ there replaced by $g_0, f_j$, respectively.
One can follow the remaining proof of \cite[Theorem 3.3]{Pe} to establish that
\[
\sum_x \bigabs{\E_{y\in[N']}f_j(x+q'q^by)}\gg_{C,d}\delta^{O_{d}(1)}N,
\]
where $q'\ll_{C,d} \delta^{-O_{d}(1)} $, $b\ll_{d} 1+ \max_{1\leq i\leq j-1}b_i \ll_{d} 1$, and 
\[
N' \gg_{C,d} (\delta/q)^{O_{d}(1)} M_{j+1}^{\deg P_j}\gg_{C,d} (\delta/q)^{O_{d}(1)} M^{\deg P_j}.
\]
 An application of the Cauchy-Schwarz inequality then yields 
\[
\delta^{O_{d}(1)} \ll_{C,d} N^{-1} \sum_{x\in \Z} \E_{y,y'\in [N]} f_j(x+q'q^b y)\overline{ f_j(x+q'q^b y')}.
\]
Making the change of variables $x+q'q^b y \to x$ then gives
\[
\delta^{O_{d}(1)} \ll_{C,d} N^{-1} \sum_{x\in \Z} f_j(x) \E_{y,y'\in [N]} \overline{ f_j(x+q'q^b (y'-y))}.
\]
The proposition follows by noting, from Lemma \ref{lipschitz}, that the function  
\[
\phi_j (x) = \E_{y,y'\in [N]} \overline{ f_j(x+q'q^b (y'-y))}
\]
 is $O_{C,d}((q/\delta)^{-O_d(1)} M^{-\deg P_j})$-Lipschitz along $q'q^b\cdot \Z$.
\end{proof}

\subsection{Correlation for all functions}\label{sec:inverse2}

In order to establish correlations for all functions, the final tool we require is a decomposition lemma due to Gowers \cite[Proposition 3.6]{Go10}. If $\|\cdot\|$ is a seminorm on an inner product space, recall that its dual seminorm $\|\cdot\|^*$ is defined by
\begin{equation}\label{dual norm}
\|f\|^{*} := \sup_{\|g\|\leq1}|\langle f,g\rangle|.
\end{equation}
Hence,
\begin{equation}\label{dual ineq}
\abs{\ang{f,g}} \leq \norm{f}^* \norm{g}.
\end{equation}

\begin{lemma}[Decomposition]\label{decomposition}
Let $\norm{\cdot}$ be a norm on the space of complex-valued functions with support on the interval $[N]$. If $\eps>0$ and $f:\Z\to\C$ is a function supported on the interval $[N]$, then there is a decomposition $f = f_{\str} + f_\unf$ such that
\[
\norm{f_\str}^* \leq \eps^{-1} \norm{f}_2 \quad \text{ and }\quad \norm{f_\unf} \leq \eps \norm{f}_2.
\]
	
\end{lemma}

%We are now ready to state the main theorem of this section, which asserts that if the counting expression $\Lambda_{P_1,\dots,P_m}^{N,M}(f_0,\dots,f_m) $ is large, then for each function $f_j$ there exists a suitable Lipschitz function that is correlated with it.
% 
%\begin{theorem}[Full inverse theorem]\label{full-inverse}
%Let $N>0$, $q\in\N$ and $P_1,\dots,P_m\in\Z[y]$ be polynomials with $(C,q)$-coefficients such that $\deg P_1< \cdots  <\deg P_m:=d$. \change{Define $M: = (N/q^{d-1})^{1/d}$. Suppose that $f_0,\dots,f_m:\Z\to\C$ are 1-bounded functions supported on the interval $[N]$.} If 
%\[
%\bigabs{ \Lambda_{P_1,\dots,P_m}^{N,M}(f_0,\dots,f_m) } \geq  \delta,
%\]
%then either $N\ll_{C,d} (q/\delta)^{O_{d}(1)}$ or, there exist positive integers $q'\ll_{C,d} \delta^{-O_{d(1)}}$ and $b' \ll_{d} 1$, as well as 1-bounded functions $\phi_i:\Z \to\C (1\leq i\leq m)$ which are $O_{C,d}((q/\delta)^{O_{d}(1)} M^{-\deg P_i})$-Lipschitz along $q'q^b\cdot \Z$ respectively, and $\phi_0:\Z\to\C$  is also  a 1-bounded,  $O_{C,d}((q/\delta)^{O_{d}(1)})M^{-\deg P_1})$-Lipschitz function along $q'q^b\cdot \Z$, such that
%\[
%\bigabs{\sum_x f_i(x)\phi_i(x)} \gg_{C,d} \delta^{O_{d}(1)}N \qquad (0\leq i\leq m).
%\]
%	
%\end{theorem}
%
%\begin{remark}
%In this theorem we assume that $\supp(f_i) =[N]$ for all $0\leq i\leq m$.	
%\end{remark}

We are now ready to prove Theorem \ref{thm:inverse} under the additional assumption that $1/C \leq q^{d-1}M^d/N \leq C$. 
Without loss of generality, we may assume for the remainder of this section that  $\deg P_1 < \cdots < \deg P_m$. Let $d_i = \deg P_i$ and $d = d_m$. We proceed by induction to prove the existence of $1$-bounded functions $\Phi_1,\cdots,\Phi_{m-1}:\Z\to\C$, where each $\Phi_i$ is $O_{C,d}((q/\delta)^{O_{d}(1)} M^{-\deg P_i})$-Lipschitz along $q'q^{b}\cdot\Z$ for some $q'\ll_{C,d}\delta^{-O_{d}(1)}$ and $b\ll_{d}1$, such that
\begin{equation}\label{eq:induct}
\bigabs{\Lambda_{P_1,\cdots,P_m}^{N,M}(f_0,\Phi_1,\cdots,\Phi_{j-1},f_j,\cdots,f_m)} \geq \delta_j
\end{equation}
for every $1 \leq j \leq m$ and some $\delta_j \gg_{C,d}\delta^{O_d(1)}$. Once this is proved, the desired correlation results for all $f_j$'s follow from Proposition \ref{partial-ii}.

The base case $j=1$ of \eqref{eq:induct} follows from Proposition \ref{partial-ii} (1) directly. For the induction step, suppose that $\Phi_1,\cdots,\Phi_{j-1}$ have been constructed for some $1 \leq j < m$ such that \eqref{eq:induct} holds. Our goal is to construct $\Phi_j$ such that \eqref{eq:induct} holds for $j+1$.

Define the semi-norm $\norm{\cdot}$ by setting
\[
\norm{f} = \sup_{g_0,g_{j+1},\cdots,g_m} N \cdot \Bigabs{\Lambda_{P_1,\cdots,P_m}^{N,M}(g_0,\Phi_1,\cdots,\Phi_{j-1},f,g_{j+1},\cdots,g_m)},
\]
where the supremum is taken over all $1$-bounded functions $g_0,g_{j+1},\cdots,g_m$ supported on $[N]$.
Applying Lemma \ref{decomposition} with $\eps=\delta_j N^{1/2}/2$ we obtain a decomposition $f_j =f_j^{\text{str}}+ f_j^{\text{unf}}$ with
\begin{align}\label{str}
\norm{f_j^{\text{str}}}^*\leq 2\delta_j^{-1} N^{-1/2} \norm{f_j}_2\leq 2\delta_j^{-1}
\end{align}
and
\begin{align}\label{unf}
\norm{f_j^{\text{unf}}}\leq (\delta_j/2)N^{1/2}\norm{f_j}_2\leq \delta_j N /2.
\end{align}
After a change of variables replacing $x$ by $x-P_j(y)$, we can write
$$
 \Lambda_{P_1,\cdots,P_m}^{N,M}(f_0,\Phi_1,\cdots,\Phi_{j-1},f_j,\cdots,f_m) = N^{-1} \langle f_j, G_j\rangle,
$$
where $G_j$ is the dual function defined by
\[
G_j(x)=\E_{y\in [M]} f_0(x-P_j(y)) \prod_{1\leq i \leq j-1} \Phi_i(x+P_i(y)-P_j(y)) \prod_{j+1 \leq i \leq m}f_i(x+P_i(y)-P_j(y)).
\]
Then it follow from \eqref{eq:induct} and the triangle inequality  that
\[
\delta_j \leq N^{-1} | \ang{f_j, G_j}| \leq N^{-1} |\ang{f_j^{\text{unf}},G_j}| + N^{-1} |\ang{f_j^{\text{str}},G_j}|.
\]
From the definition of the norm $\norm{\cdot}$ and the inequality (\ref{unf}) one can deduce that
 \[
 |\ang{f_j^{\text{unf}},G_j}|\leq \norm{f_j^{\text{unf}}}\leq \delta_j N/2.
 \] 
Hence, it follows from (\ref{dual ineq}), (\ref{str}) and the above two inequalites that
\[
\delta_j/2 \leq N^{-1} |\ang{f_j^{\text{str}},G_j}| \leq N^{-1} \norm{f_j^{\text{str}}}^* \norm{G_j} \leq 2\delta_j^{-1} N^{-1} \norm{G_j}.
\]
Therefore, $\norm{G_j}\gg \delta_j^2 N$, which implies that there exist 1-bounded functions $g_0,g_{j+1},\cdots,g_m:\Z\to\C$ supported on the interval $[N]$ such that
\[
\bigabs{\Lambda_{P_1,\dots, P_m}^{N,M}(g_0,\Phi_1,\cdots,\Phi_{j-1},G_j,g_{j+1},\cdots,g_m} \gg \delta_j^2 \gg_{C,d}\delta^{O_d(1)}.
\]
By Proposition \ref{partial-ii}, there exists a $1$-bounded function $\Phi_j:\Z\to\C$ which is $O_{C,d}((q/\delta)^{O_{d}(1)} M^{-\deg P_j})$-Lipschitz along $q_jq^{b_j}\cdot Z$ for some $q_j\ll_{C,d}\delta^{-O_{d}(1)}$ and $b_j\ll_{d}1$, such that
\[
\Bigabs{\sum_x G_j(x) \Phi_j(x)} \gg_{C,d} \delta^{O_{d}(1)}N.
\]
We may enlarge $q',b$ if necessary to ensure that all of $\Phi_1,\cdots,\Phi_j$ are Lipschitz functions along $q'q^b \cdot \Z$, while maintaining the bounds $q' \ll_{C,d}\delta^{-O_d(1)}$ and $b\ll_d1$.
Expanding the dual function $G_j$, we thus obtain that
\begin{align}\label{bias-f2}
\Bigabs{\Lambda_{P_1,\cdots,P_m}^{N,M}(f_0,\Phi_1,\cdots,\Phi_{j-1},\Phi_j,f_{j+1},\cdots,f_m)} \gg_{C,d} \delta^{O_{d}(1)}.
\end{align}
This completes the induction step, thereby proving Theorem \ref{thm:inverse} in the case when $1/C \leq q^{d-1}M^d/N \leq C$.

\begin{remark}
In the above argument, we used Lemma \ref{decomposition} based on the Hahn-Banach theorem to deduce intermediate conclusions for the dual functions $G_j$, which we then use to deduce conclusions for $f_{j+1}$. A modern way of executing this maneuver (in a much more sophisticated setting) was recently introduced by Manners \cite{Manners} and referred to as ``stashing".
\end{remark}

\subsection{Theorem \ref{thm:inverse} for general $M$}\label{sec:inverse3}

Finally, we now address the general case of Theorem \ref{thm:inverse} by reducing to the case considered in Section \ref{sec:inverse2} after dividing $[N]$ into subintervals of appropriate lengths.

Since each $P_j$ has $(C,q)$-coefficients and $\deg P_j\leq d$, we have for $1 \leq y \leq M$,
\begin{equation}\label{eq:sec2.3-1}
|P_j(y)| \leq d C^2 q^{d-1} |y|^d \leq d C^2 q^{d-1}M^d.
\end{equation}
Set $N_0 =  dC^2 q^{d-1}M^d$. We may assume that $N_0 \leq N$, since otherwise $N/(dC^2) \leq q^{d-1}M^d \leq N$ and the conclusion follows from the previous case (after replacing $C$ by $dC^2$). Now divide $[N]$ into a collection $\mathcal{I}$ of $\asymp N/N_0$ subintervals, such that each interval $I \in \mathcal{I}$ has length between $N_0$ and $2N_0$.

For the rest of the proof, fix some $1 \leq i \leq m$, and  we aim at deducing the correlation condition for $f_i$. Writing $f_i\vert_I$ for the restriction of $f_i$ on $I$, we have
\begin{equation}\label{eq:sec2.3-2}
\bigabs{\Lambda_{P_1,\dots,P_m}^{N,M}(f_0, \cdots,f_{i-1}, f_i\vert_{I}, f_{i+1},\dots,f_m)} \gg \delta \cdot \frac{N_0}{N}
\end{equation}
for each interval $I$ in a subcollection $\mathcal{I}'$ of $\mathcal{I}$ with $|\mathcal{I}'| \gg N/N_0$. Note that if $x+P_i(y) \in I$, then from \eqref{eq:sec2.3-1} we have $x \in J$ and $x + P_j(y) \in J$ for each $1 \leq j \leq m$, where $J = I + [-2N_0, 2N_0]$ is an interval of length at most $6N_0$. Hence we may assume that all of the functions $f_0,\dots,f_{i-1},f_i\vert_{I},f_{i+1},\dots,f_m$ appearing in \eqref{eq:sec2.3-2} are supported on $J$. By defining $\tilde{f}_j$ to be a suitable shift of $f_j$ for $j \neq i$ and defining $\tilde{f}_i$ to be a suitable shift of $f_i\vert_I$, we may ensure that each $\tilde{f}_j$ is supported on $[1, 6N_0]$ and obtain
\[
\bigabs{\Lambda_{P_1,\dots,P_m}^{6N_0,M}(\tilde{f}_0, \tilde{f}_1,\cdots,\tilde{f}_m)} \gg \delta.
\]
Since $q^{d-1}M^d\asymp_C N_0$, it follows from the conclusion in Section \ref{sec:inverse2} that there exist positive integers $q_I'\ll_{C,d} \delta^{-O(1)}$ and $b_I\ll 1$, as well as a $1$-bounded function $\phi_I: \Z\rightarrow\C$ which is $O_{C,d}((q/\delta)^{O(1)}M^{-d_i})$-Lipschitz along $q_I'q^{b_I}\cdot \Z$, such that
\begin{equation}\label{eq:phi_I}
\bigabs{\sum_{x \in \Z} f_i\vert_I(x) \phi_I(x)} \gg_{C,d} \delta^{O_d(1)}N_0.
\end{equation}
By the pigeonhole principle, there exists a subset $\mathcal I'' \subset \mathcal I'$ satisfying $|\mathcal{I}''| \gg_{C,d} \delta^{O_d(1)} N/N_0$, such that $q_I' = q'$  and $b_I = b$ are independent of $I$ for $I \in \mathcal{I}''$. 

For $I \in \CI''$, we may change the values of $\phi_I$ near the endpoints of $I$ appropriately to obtain a modified function $\phi_I': \Z\rightarrow\C$ which is supported on $I$ and takes the value $0$ at the endpoints of $I$, while ensuring that $\phi_I'$ remains $O_{C,d}((q/\delta)^{O(1)}M^{-d_i})$-Lipschitz along $q'q^{b}\cdot \Z$ and that \eqref{eq:phi_I} still holds with $\phi_I$ replaced by $\phi_I'$. More specifically, let
$$
\eta = \frac{1}{10N_0} \bigabs{\sum_{x \in \Z} f_i\vert_I(x) \phi_I(x)}  \gg_{C,d} \delta^{O_d(1)}.
$$
If $I = [u, v]$ then we define $\phi_I'$ by $\phi_I'(x)=0$ for $x \leq u$ and for $x \geq v$, $\phi_I'(x) = \phi_I(x)$ for $x \in [u+\eta N_0, v-\eta N_0]$, and requiring that $\phi_I'$ is linear on $[u, u+\eta N_0]$ and on $[v-\eta N_0, v]$. By linearity, $\phi_I'$ is $O(\eta^{-1}N_0^{-1})$ Lipschitz on  $[u, u+\eta N_0]$ and $[v-\eta N_0, v]$, and thus the desired Lipschitz condition on $\phi_I'$ is satisfied. Moreover, replacing $\phi_I$ by $\phi_I'$ in \eqref{eq:phi_I} does not change the summand unless $x \in [u, u+\eta N_0]$ or $x \in [v-\eta N_0, v]$. Since there are at most $2\eta N_0$ such values of $x$, it follows that \eqref{eq:phi_I} remains valid for $\phi_I'$ by our choice of $\eta$.

Next, by replacing $\phi_I'$ by $z_I\phi_I'$ for some complex number $z_I$ with $|z_I|=1$, we may assume that 
$$
\sum_{x \in \Z} f_i\vert_I(x) \phi_I'(x) \gg_{C,d} \delta^{O(1)}N_0.
$$
Finally, define $\phi:\Z\rightarrow\C$ by gluing the functions $\phi_I'$ together for $I \in \mathcal{I}''$; i.e. set $\phi(x) = \phi_I'(x)$ if $x \in I$ for some $I \in \CI''$ and set $\phi(x)=0$ otherwise. This establishes the desired conclusion, thereby completing the proof of Theorem \ref{thm:inverse}.

%%%%%%%%%%%%%%%%%%%%%%%%%%%%%%%%%%%%%%%%%%%%%%%%%%%%%%%%%%%%%%%%%%%%%%%%%%%%%%%%%%%%%%%%%%%%%%%%%%%%%%%%%%%%%%%%%%%%%%%%%%%%%%

\section{Local factors}\label{sec:factors}

In this section, we adopt the notation for factors from \cite{GT-4-aps}, originally derived from ergodic theory; see also \cite[Section 4]{PP2}.

\begin{definition}[Factor]
  We define a \emph{factor} $\mathcal{B}$ of $[N]$ to be a partition of $[N]$, so that $[N] = \sqcup_{B \in \mathcal{B}} B$.   
  A factor $\mathcal{B}'$ \emph{refines} $\mathcal{B}$ if every element of $\mathcal{B}$ is a union of elements of $\mathcal{B}'$.  
  The \emph{join} $\mathcal{B}_1\vee\dots\vee\mathcal{B}_d$ of factors $\mathcal{B}_1, \dots, \mathcal{B}_d$ is the factor formed by taking the $d$-fold intersections of the elements of $\mathcal{B}_1$, \dots, $\mathcal{B}_d$, that is,
  \[
\mathcal{B}_1\vee\dots\vee\mathcal{B}_d:=\{B_1\cap\dots\cap B_d:B_i\in\mathcal{B}_i\text{ for }i=1,\dots,d\}.
  \]
\end{definition}
\begin{definition}[Measurability, projection]
Given a factor $\mathcal{B}$ of $[N]$, we say that a function $f : [N] \to \C$ is \emph{$\mathcal{B}$-measurable} if it is constant on the elements of $\mathcal{B}$.  
Define the \emph{projection} of any function $f : [N] \to \C$ onto $\mathcal{B}$  by
\begin{equation}\label{conditional expectation}
\Pi_\mathcal{B} f(x)  = \E_{y \in B(x)} f(y),
\end{equation}
where $B(x)$ is the unique atom of $\mathcal{B}$ containing $x$. 

\end{definition}

Notice that $\Pi_\mathcal{B} f  $ is $\mathcal{B}$-measurable and corresponds to the conditional expectation of $\mathcal{B}$ with respect to the $\sigma$-algebra generated by the elements of $\CB$.
The following lemma illustrates why we work with projections onto factors:  the $l^2$-norm of such projections is non-decreasing under refinement of the factor. For this reason, we refer to $\norm{\Pi_\CB f}_2^2$ as the \emph{energy}.

\begin{lemma}[Pythagoras theorem for projections]\label{pythogoras}
Let $\CB,\CB'$ be factors of $[N]$ such that $\CB'$ refines $\CB$. Let $f: [N]\rightarrow \C$ be a function. Then
$$
\bignorm{\Pi_{\CB'}f}_2^2 = \bignorm{\Pi_{\CB}f}_2^2 + \bignorm{\Pi_{\CB'}f -\Pi_\CB f}_2^2.
$$
\end{lemma}

\begin{proof}
See \cite[Lemma 4.3]{PP2}.
\end{proof}

In this paper we will exclusively work with factors whose atoms are arithmetic progressions.

\begin{definition}[Local factor]\label{local factor def}
Let $\CB$ be a factor of $[N]$. Let $q, M \geq 1$.  We say that $\CB$ is a  \emph{local factor} of resolution $M$ and modulus $q$ if every atom of $\CB$ is an arithmetic progression of step $q$ and length in $[M, 2M]$.\end{definition}

For example, the trivial factor consisting of only one atom is a local factor of resolution $N$ and modulus $1$. Our definition of local factors is motivated by, but slightly different from, the definition in \cite{PP2}, where the atoms are required to have length exactly $M$. This flexibility leads to the following simple lemma which turns out to be rather convenient for our arguments.  

\begin{lemma}[Existence of local factors]\label{lem:local-factor}
Let $\CB, \CB''$ be  local factors of $[N]$ of modulus $q$ and resolution $M, M''$, respectively, and $100M''\leq M$. Suppose that $\CB''$ refines $\CB$. Then for any integer $M' \in [10M'', M/10]$, there exists a  local factor $\CB'$ of $[N]$ of resolution $M'$ and modulus $q$, such that $\CB'$ refines $\CB$ and $\CB''$ refines $\CB'$.
\end{lemma}

\begin{proof}
Let $P$ be an arbitrary atom of $\CB$, so that $P$ is an arithmetic progression of step $q$ and length in $[M, 2M]$. Since $\CB''$ refines $\CB$, we have a partition 
$$ P = P_1 \sqcup P_2 \sqcup \cdots \sqcup P_r, $$
where each $P_i$ is an atom of $\CB''$ which is an arithmetic progression of step $q$ and length in $[M'', 2M'']$. Without loss of generality, we may assume that all elements of $P_i$ are smaller than all elements of $P_j$ whenever $i<j$. Let $x_i = |P_i|/M'$ so that $x_i \leq [y, 2y]$ where $y = M''/M' \leq 1/10$. Then
$$ s := x_1 + x_2 + \cdots + x_r = \frac{|P|}{M'} \geq \frac{M}{M'} \geq 10. $$ 
We will define a sequence
$$ 0 = i_0 < i_1 < \cdots < i_{k-1} < i_k = r, $$
such that 
$$ x_{i_{j-1}+1} + \cdots + x_{i_j} \in [1, 2] $$
for each $1 \leq j \leq k$. Once this construction is completed, we can partition $P$ into $k$ arithmetic progressions $P_{i_{j-1}+1} \cup \cdots \cup P_{i_j}$ ($1 \leq j \leq k$), each of which has step $q$ and length in $[M', 2M']$. Performing this procedure for each atom of $\CB$, we obtain the desired refinement $\CB'$ of $\CB$.

To construct the sequence $\{i_j\}$, first choose a positive integer $k$ such that $s/k \in [1.4, 1.6]$. The existence of such $k$ easily follows from $s \geq 10$. Now, for each $1 \leq j \leq k$, define $i_j$ to be the smallest index with the property that
$$ x_1 + x_2 + \cdots + x_{i_j} \geq \frac{js}{k}. $$
Clearly we must have $i_k = r$ and the upper bound
$$ x_1 + x_2 + \cdots + x_{i_j} \leq x_1 + \cdots + x_{i_j-1} + 2y \leq \frac{js}{k} + 0.2. $$
It follows that
$$ x_{i_{j-1}+1} + \cdots + x_{i_j} \leq \frac{js}{k}+0.2 - \frac{(j-1)s}{k} = \frac{s}{k} + 0.2 \leq 1.8 $$
and
$$ x_{i_{j-1}+1} + \cdots + x_{i_j} \geq \frac{js}{k} - \left(\frac{(j-1)s}{k}+0.2\right) = \frac{s}{k}-0.2 \geq 1.2. $$
This completes the proof.
\end{proof}

In our proof of the popular difference result (Theorem \ref{thm:popular}), we need to work with a chain of local factors where each factor in the chain refines the next one.

\begin{definition}[Local factor chain]\label{def:local-factor-chain}
Let $q, M_1,\dots,M_m$ be positive integers.
Let $\CB_1,\dots,\CB_m$ be local factors of $[N]$ of modulus $q$ and resolution $M_1,\dots,M_m$, respectively. We say that $(\CB_1,\dots,\CB_m)$ is \emph{a local factor chain} of resolution $(M_1,\dots,M_m)$ and modulus $q$, if $\CB_i$ is a refinement of $\CB_{i+1}$ for each $1 \leq i < m$. 

Let $(\CB_1,\dots,\CB_m)$ and $(\CB_1',\dots,\CB_m')$ be two local factor chains. We say that $(\CB_1',\dots,\CB_m')$ \emph{refines} $(\CB_1,\dots,\CB_m)$ if $\CB_i'$ refines $\CB_i$ for each $1 \leq i \leq m$.
\end{definition}

To visualize the relationships between factors, we use the notation $\CB \longrightarrow \CB'$ to indicate that $\CB$ refines $\CB'$. Thus a local factor chain $(\CB_1,\cdots,\CB_m)$ can be visualized as:
$$
\xymatrix{
\CB_1 \ar[r] & \CB_2 \ar[r] & \cdots \ar[r] & \CB_m 
}
$$
and the relationship that $(\CB_1',\cdots,\CB_m')$ refines $(\CB_1,\cdots,\CB_m)$ can be visualized as:
$$
\xymatrix{
\CB_1 \ar[r] & \CB_2 \ar[r] & \cdots \ar[r] & \CB_m \\
\CB_1' \ar[u] \ar[r] & \CB_2' \ar[u]\ar[r] & \cdots \ar[r] & \CB_m' \ar[u]
}
$$

\begin{lemma}[Existence of local factor chains]\label{lem:local-factor-chain}
Let $(\CB_1,\dots,\CB_m)$ be a local factor chain of $[N]$ of resolution $(M_1,\dots,M_m)$ and modulus $q$. Let $q', M_1',\dots,M_m'$ be positive integers such that $q'$ is divisible by $q$, and
$$ 10M_{i-1} \leq \frac{q'}{q}M_i' \leq \frac{M_i}{10} $$
for each $1 \leq i \leq m$, with the convenience that $M_0=1$.
Then there exists a local factor chain $(\CB_1',\dots,\CB_m')$ of resolution $(M_1',\dots,M_m')$ and modulus $q'$ such that $(\CB_1',\dots,\CB_m')$  refines $(\CB_1,\dots,\CB_m)$.
\end{lemma}

\begin{proof}
First we prove the lemma in the special case $q'=q$.
For each $1 \leq i \leq m$, we will apply Lemma \ref{lem:local-factor} with $\CB = \CB_{i}$, $\CB'' = \CB_{i-1}$ (taking $\CB_0$ to be the trivial factor where each atom is a singleton), $M = M_i$, $M'' = M_{i-1}$, and $M' = M_i'$. The hypothesis $10M'' \leq M' \leq M/10$ is satisfied by our assumption. Hence we obtain a local factor $\CB_i'$ of resolution $M_i'$ and modulus $q$, such that $\CB_i'$ refines $\CB_i$ and $\CB_{i-1}$ refines $\CB_i'$. The following diagram illustrates the construction in the case $m=3$, where the dotted arrows represent refinement relations from applying Lemma \ref{lem:local-factor}:
\[
\xymatrix{
\CB_0 \ar[r] \ar @{.>}[dr] & \CB_1 \ar[r]  \ar @{.>}[dr] & \CB_2 \ar[r] \ar @{.>}[dr] & \CB_3  \\
 & \CB_1'\ar @{.>}[u]^{Lemma \ref{lem:local-factor}} \ar @2{=>}[r]  & \CB_2'\ar @{.>}[u]^{Lemma \ref{lem:local-factor}} \ar @2{=>}[r]  & \CB_3'\ar @{.>}[u]^{Lemma \ref{lem:local-factor}}
}
\]
Since $\CB_{i-1}'$ refines $\CB_{i-1}$ and $\CB_{i-1}$ refines $\CB_i'$, it follows that $\CB_{i-1}'$ refines $\CB_i'$ for all $2\leq i\leq m$, as illustrated in the diagram above by the double arrows. Therefore, $(\CB_1',\dots,\CB_m')$ is a local factor chain which refines $(\CB_1,\cdots,\CB_m)$, as desired.

Now we treat general $q'$. By the above procedure, we first obtain a local factor chain $(\CB_1'',\dots,\CB_m'')$ of resolution $(q'M_1'/q, \cdots, q'M_m'/q)$ and modulus $q$ which refines $(\CB_1,\dots,\CB_m)$. Next, for each $1 \leq i \leq m$, divide each atom of $\CB_i''$ into residue classes modulo $q'$ to form a refinement $\CB_i'$ of $\CB_i''$. Since the atoms of $\CB_i''$ have length in $[q'M_i'/q, 2q'M_i'/q]$  and have step $q$, the atoms of $\CB_i'$ have length in $[M_i', 2M_i']$. Hence $(\CB_1',\dots,\CB_m')$ is a local factor chain of resolution $(M_1',\dots,M_m')$ and modulus $q'$ which refines $(\CB_1'',\cdots,\CB_m'')$. This completes the proof.
\end{proof}

%%%%%%%%%%%%%%%%%%%%%%%%%%%%%%%%%%%%%%%%%%%%%%%%%%%%%%%%%%%%%%%%%%%%%%%%%%%%%%%%%%%%%%%%%%%%%%%%%%%%%%%%%%%%%%%%%%%%%%%%%%%%%%
\section{Polylogarithmic bound in the density result}\label{sec:density-increment}

In this section, we prove Theorem \ref{thm:density} using the improved density increment strategy developed by Heath-Brown \cite{HB-3aps} and Szemer\'edi \cite{sze-3aps}.  For convenience, we follow the presentation in Green--Tao \cite{GT-4-aps}.  

We begin by briefly outlining the  idea. Assume that $A\subseteq[N]$ is a subset of density $\alpha$ that contains no nontrivial polynomial configurations. Since Lipschitz functions are nearly constant on long arithmetic progressions, Peluse \cite{Pe} deduced from the inverse theorem (Theorem \ref{thm:inverse}) and pigeonhole principle that that there exists a progression $P$ on which $A$ exhibits a density increment of size $\alpha^{O(1)}$. It turns out that this approach typically requires $O(\alpha^{-O(1)})$ iterations to reach a contradiction. The strategy adopted in this section improves on this by collecting  several Lipschitz functions together to construct a new one, and then decomposing $[N]$ into progressions on which this composite function is nearly constant. This refinement allows us to obtain a density increment of size $c\alpha$ on one of the progressions, thereby reducing the number of iterations to just $O(\log 1/\alpha)$. While the progressions used at each step are shorter than those in Peluse's approach, it was observed in \cite{GT-4-aps} that the progression length plays a less critical role in the overall argument.

Let us clarify the above argument. We first construct a nonnegative structural function $g\geq 0$ such that it has the same mean value as  $1_A$ (say, $\alpha$), and the corresponding counting expressions weighted by $g$ and $1_A$ are comparable. This construction is formalized in Lemma \ref{weak-regularity} below. We then carry out the density-increment argument with respect to $g$. Lemma \ref{density-increment}, in essence, asserts that the density of $g$ increases by a factor of at least $(1+c)\alpha$ on some arithmetic progression.

 \begin{lemma}[Weak regularity lemma-I]\label{weak-regularity}
  Let $M, N,q$ be positive integers and let $P_1,\dots,P_m\in\Z[y]$ be polynomials with $(C,q)$-coefficients such that $\deg P_1< \cdots  <\deg P_m$. Suppose that $f:\Z\to\C $ is a 1-bounded function  supported on the interval $[N]$. Let $d = \deg P_m$, and assume that $M \leq (N/q^{d-1})^{1/d}$. For $\delta \in (0,1/2)$, one of the following two statements holds:

\begin{enumerate}
\item $M\ll_{C,d}(q/\delta)^{O_{d}(1)}$.

\item There exist positive integers $Q, b, M_1,\dots,M_m$ with
  $$
  Q \leq \exp(O_{C,d}(\delta^{-O_d(1)})), \ \ b\ll_d1, \ \ M_i\gg_{C,d} Q^{-1}(\delta/q)^{O_d(1)}M^{\deg P_i},
  $$
and local factors $\CB_1,\dots,\CB_m$  on $[N]$ of resolution $M_1,\dots,M_m$, respectively, and modulus $Qq^b$, such that
  \[
  \bigabs{ \Lambda_{P_1,\dots,P_m}^{N,M} (f,\dots, f) - \Lambda_{P_1,\dots,P_m}^{N,M} (f, \Pi_{\CB_1}f,\dots, \Pi_{\CB_m}f)}\leq \delta.
  \] 
\end{enumerate}
    \end{lemma}
 
%\change{It is not necessary for $(\CB_1,\dots,\CB_m)$ to form a local factor chain in the application of this lemma (i.e. $\CB_i$ refines $\CB_{i+1}$). This allows us to perform the iteration steps for $i=1,\dots,m$, independently. We only need to take $q=\lcm\set{q_1,\dots,q_m}$ where $q_i$ is the modulus of $\CB_i$ obtained via the energy increment argument for each  $i$. I hope this observation helps simplify the proof.}

\begin{proof}
Let $b' = b'(d)$ be a positive integer chosen sufficiently large  in terms of $d$, and let $C' = C'(C,d)$ be a constant taken large enough in terms of $C,d$. Set $Q$ to be the least common multiple of all positive integers at most $C'\delta^{-b'}$, and set $M_i = C'^{-1}Q^{-1}(\delta/q)^{b'} M^{\deg P_i}$. The desired bounds for $Q$ and $M_i$ are then satisfied.

Let $\CB_1,\dots,\CB_m$ be local factors on $[N]$ of resolution $M_1,\dots,M_m$, respectively, and modulus $Qq^b$. Decompose $\Lambda_{P_1,\dots,P_m}^{N,M}(f,\dots,f)$ into $2^m$ terms of the form
$$
\Lambda_{P_1,\dots,P_m}^{N,M}(f, g_1,\dots,g_m), 
$$
where each $g_i \in \{\Pi_{\CB_i}f, f - \Pi_{\CB_i}f\}$. It suffices to show that if $g_i = f - \Pi_{\CB_i}f$ for some $1 \leq i \leq m$ then
$$
|\Lambda_{P_1,\dots,P_m}^{N,M}(f, g_1,\dots,g_m)| \leq \frac{\delta}{2^m}.
$$
Suppose, on the contrary, that this is not the case. Then by Theorem \ref{thm:inverse}, either $M\ll_{C,d}(q/\delta)^{O_{d}(1)}$ in which case we are done, or there exist positive integers $q' \ll_{C,d}\delta^{-O_d(1)}$ and $b \ll_d1$, as well as a $1$-bounded function $\phi_i:\Z\rightarrow\C$ which is $L$-Lipschitz along $q'q^b\cdot \Z$ for some $L \ll_{C,d}(q/\delta)^{O_d(1)}M^{-\deg P_i}$, such that
$$
\bigabs{ \sum_{x \in \Z} g_i(x)\phi_i(x) } \geq \eta N
$$
for some $\eta \gg_{C,d} \delta^{O_d(1)}$.
We may ensure that $b$ is a constant depending only on $d$ and that $q' \mid Q$ by choosing $b', C'$ large enough. If $x,y$ lie in the same atom $P$ of $\CB_i$, then $q'q^b \mid x-y$ and $|x-y| \leq 2Qq^b M_i$. Hence
$$
|\phi_i(x) - \phi_i(y)| \leq L|x-y| \leq  2L Q q^b M_i,
$$
which we may ensure to be at most $\eta/2$ by choosing $b',C'$ large enough. Since $P$ is an atom of $\CB_i$,  the average of $g_i=f-\Pi_{\CB_i}f$ on $P$ is $0$. It follows that
$$
\bigabs{\sum_{x \in P} g_i(x)\phi_i(x)} \leq \frac{1}{2}\eta |P|
$$
for each atom $P$ of $\CB_i$. This leads to a contradiction after summing over all atoms and concludes the proof.
\end{proof}

\begin{lemma}[Density increment]\label{density-increment}

 Let $M, N,q$ be positive integers and $P_1,\dots,P_m\in\Z[y]$ be polynomials with $(C,q)$-coefficients such that $\deg P_1< \cdots  <\deg P_m$.  Let $d = \deg P_m$ and $M = (N/q^{d-1})^{1/d}$. Suppose that $ A \subseteq [N]$ has density $\alpha := | A|/N$ and contains no nontrivial progression of the form $x,x+P_1(y),\dots, x+P_m(y)$. Then either $N\ll_{C,d} (q/\alpha)^{O_d (1)}$ or there exist integers $q'\leq \exp \bigbrac{ O_{C,d} (\alpha^{-O_d (1)})} $, $b \ll_d 1$ and an arithmetic progression $P\subseteq[N]$ of modulo $q'q^b$ and of length $\gg_{C,d}(\alpha/q)^{O_d(1)}M^{\deg P_1} /q'$  such that
\[
| A \cap P| \geq \alpha (1+c) |P|
\]
for some constant $c=c(C,d) > 0$.	
\end{lemma}
 
Let us state the following  $l^1$-control lemma which will be used in proving Lemma \ref{density-increment}. This lemma is an extension of \cite[Lemma 5.1]{PP2}, and its proof is elementary.

\begin{lemma}[$l^1$-control]\label{l1-control}
For any functions $f_0,\dots,f_m:[N]\to\C$ we have
\[
\bigabs{ \Lambda_{P_1,\dots,P_m}^{N,M} (f_0,\dots, f_m)}\leq N^{-1} \norm{f_i}_1 \prod_{j\neq i} \norm{f_j}_\infty.
\]	
\end{lemma}

\begin{proof}[Proof of Lemma \ref{density-increment}]
By the assumption on $A$, we have $\Lambda_{P_1,\dots,P_m}^{N,M}(1_A,\dots,1_A)=0$.  Besides, it is obvious that
$$
\Lambda_{P_1,\dots,P_m}^{N,M}(1_A,\alpha 1_{[N]},\dots,\alpha 1_{[N]}) = \frac{\alpha^m}{NM} \sum_{x \in \Z} 1_A(x) \sum_{y \in [M]}  1_{[N]}(x+P_1(y)) \cdots 1_{[N]}(x+P_m(y)).
$$
Let $\eta > 0$ be a constant sufficiently small in terms of $C,d$ so that $|P_i(y)| \leq N/10$  whenever $1 \leq y \leq \eta M$ and $1\leq i\leq m$. 
It follows that
\begin{align}\label{trivial-lower-bound}
	&\Lambda_{P_1,\dots,P_m}^{N,M}(1_A,\alpha 1_{[N]},\dots,\alpha 1_{[N]})\nonumber\\
	&\geq \frac{\alpha^m}{NM} \sum_{x\in [N/3,2N/3]} 1_A(x) \sum_{y\leq \eta M} 1_{[N]} (x+P_1(y))\dots 1_{[N]} (x+P_m(y))\nonumber\\
	& \geq \frac{\alpha^m}{NM} \Big|A \cap [N/3, 2N/3]\Big|\cdot \eta M.
\end{align}
If $|A \cap [N/3, 2N/3]| \leq \alpha N/10$, it follows from the pigeonhole principle that either $|A \cap [1,N/3| \geq \alpha N/5$, or $|A \cap [ 2N/3, N]| \geq \alpha N/5$, and the conclusion  follows by taking either $P = [1, N/3]$ or $P = [2N/3, N]$. Hence we may assume that $|A \cap [N/3, 2N/3]| \geq \alpha N/10$ in the following. In light of (\ref{trivial-lower-bound}) one has
$$
\Lambda_{P_1,\dots,P_m}^{N,M}(1_A,\alpha 1_{[N]},\dots,\alpha 1_{[N]}) \geq 
\frac{\eta}{10}\alpha^{m+1}.
$$

On the other hand, by the weak regularity lemma (Lemma \ref{weak-regularity}) applied to $f=1_A$ and $\delta = \eta \alpha^{m+1}/20$, we can find local factors $\CB_1,\dots,\CB_m$ on $[N]$ of resolution $M_1,\dots,M_m$, respectively, and modulus  $q'q^b$ for some $q'\leq \exp(O_{C,d}(\alpha^{-O_d(1)}))$, $b\ll_d1$ and $M_i\gg_{C,d} (\alpha/q)^{O_{d}(1)} M^{\deg P_i}/q'$ such that
  \[
  \bigabs{ \Lambda_{P_1,\dots,P_m}^{N,M} (1_A, \Pi_{\CB_1}1_A,\dots, \Pi_{\CB_m}1_A)}\leq \frac{\eta \alpha^{m+1}}{20}.
  \] 
We may decompose $\Lambda_{P_1,\dots,P_m}^{N,M}(1_A,\alpha 1_{[N]},\dots,\alpha 1_{[N]})$ into $2^m$ terms, each of which takes the form
$$
\Lambda_{P_1,\dots,P_m}^{N,M}(1_A,g_1,\dots,g_m),
$$
where each $g_j \in \{\Pi_{\CB_j}1_A, \alpha 1_{[N]} - \Pi_{\CB_j}1_A\}$. It follows that 
\begin{align}\label{contradiction-sec.5}
\bigabs{ \Lambda_{P_1,\dots,P_m}^{N,M} (1_A, g_1,\dots,g_m)} \geq \frac{\eta \alpha^{m+1}}{20\cdot 2^m}
\end{align}
for some choice of $g_1,\dots,g_m$, at least one of which (say $g_i$) is $g_i = \alpha 1_{[N]} - \Pi_{\CB_i}1_A$.

Assume, for the sake of contradiction, that the desired density increment does not occur. Then the density of $A$ on each atom of each $\CB_j$ is at most $\alpha(1+c)$, so that $\|g_j\|_{\infty} \leq \alpha(1+c) \leq 2\alpha$ for every $j$. Moreover, since $g_i = \alpha 1_{[N]} - \Pi_{\CB_i}1_A$ takes constant values on atoms $P$ of $\CB_i$, we have
$$
\|g_i\|_1 =\sum_{P\in \CB_i} \Bigabs{ \sum_{x\in P}\Pi_{\CB_i} 1_A(x)- \alpha|P|}= \sum_{P\in \CB_i} \bigabs{|A\cap P| - \alpha |P|}. 
$$
On the other hand, we observe that the assumption  $\E_{x\in [N]}1_A(x) =\alpha$ implies that $\sum_{x\in [N]} g_i(x)=0$. From these two formulas, we then deduce that
$$
\|g_i\|_1 =  2 \sum_{P\in \CB_i} \max\bigset{|A \cap P| - \alpha |P|, 0} \leq 2c\alpha N.
$$
Thus, it follows from Lemma \ref{l1-control} that
$$
\bigabs{ \Lambda_{P_1,\dots,P_m}^{N,M} (1_A, g_1,\dots,g_m)} \leq N^{-1} \|g_i\|_1 \prod_{j \neq i} \|g_i\|_{\infty} \leq 2^{m+1} c \alpha^{m+1}.
$$
This contradicts the inequality (\ref{contradiction-sec.5}) if $c$ is sufficiently small. 
\end{proof}

\begin{proof}[Proof of Theorem \ref{thm:density}]

For $k \geq 0$, we iteratively construct positive integers $N_k, q_k$, polynomials $P_i^{(k)}$ for $1 \leq i \leq m$, each with $(C,q_k)$-coefficients and satisfying $\deg P_i^{(k)} = \deg P_i$, and subsets $A_k \subset [N_k]$ of density $\alpha_k := |A_k|/N_k$ which contains no nontrivial progressions of the form $x, x+P_1^{(k)}(y),\dots,x+P_m^{(k)}(y)$. 

For convenience, we write $d=\deg P_m$ and $d_i=\deg P_i$ for each $i$. We begin the iteration with $N_0 = N$, $q_0=1$, $P_i^{(0)} = P_i$ for $1 \leq i \leq m$, $A_0 = A$. At step $k$, either we have $N_k \ll_{C,d} (q_k/\alpha_k)^{O_d(1)}$ in which case we terminate the iteration process, or, by Lemma \ref{density-increment}, there exist $q' \leq \exp(O_{C,d}(\alpha_k^{-O_d(1)}))$, $b \ll_d 1$ and an arithmetic progression $P \subset [N_k]$ of modulo $q = q'q_k^b$ and of length $\gg_{C,d} q'^{-1}(\alpha_k/q_k)^{O_d(1)} N_k^{d_1/d}$ such that
\begin{align}\label{increment-k}
|A_k \cap P| \geq \alpha_k(1+c)|P|,
\end{align}
for some constant $c=c(C,d) > 0$. Let $N_{k+1} = |P|$ and write $P = q \cdot [N_{k+1}] + m$ for some $m \in \Z$. Define
$$
A_{k+1} := \{x \in [N_{k+1}]: qx+m \in A_k\}.
$$
Then inequality (\ref{increment-k}) yields that $\frac{|A_{k+1}|}{N_{k+1}} :=\alpha_{k+1} \geq \alpha_k(1+c)$. Since $A_k$ contains no nontrivial progressions of the form
$$
qx+m, qx+m+P_1^{(k)}(qy), \cdots, qx+m + P_m^{(k)}(qy),
$$
it follows that $A_{k+1}$ contains no nontrivial progressions of the form
$$
x, x+P_1^{(k+1)}(y), \cdots, x + P_m^{(k+1)}(y),
$$
where $P_i^{(k+1)}(y) := P_i^{(k)}(qy)/q$ has $(C,q_{k+1})$-coefficients with $q_{k+1}=q_kq = q_k^{b+1}q'$.

Suppose that the process above terminates at step $K$. Since
$$ 
1 \geq \alpha_K \geq (1+c)^K\alpha,
$$
we have
$$
K \leq \frac{\log (1/\alpha)}{\log (1+c)} \ll_{C,d} \log \frac{1}{\alpha}.
$$
From $q_{k+1} = q_k^{b+1}q'$ it follows that
$$
q_K \leq (q')^{1+(b+1)+\cdots+(b+1)^{K-1}} \leq (q')^{(b+1)^K} \leq \exp(-O_{C,d}(\alpha^{O_d(1)})).
$$
Since for any $1\leq k\leq K$ we have
$$
N_{k+1} \gg q_k^{-O_d(1)} \exp(-O_{C,d} (\alpha_k^{-O_d(1)}))N_k^{d_1/d} \gg \exp(-O_{C,d} (\alpha_k^{-O_d(1)}))N_k^{d_1/(2d)},
$$
it follows that
$$
N_K \gg  \exp(-O_{C,d}(\alpha^{-O_d(1)})) N^{d_1/(2d)}.
$$
Since $N_K \ll (q_K/\alpha_K)^{O_d(1)}\ll \exp(O_{C,d}(\alpha^{-O_d(1)}))$, it follows that
$$
N^{d_1/(2d)} \ll \exp(O_{C,d}(\alpha^{-O_d(1)})),
$$
which implies that $\alpha \ll (\log N)^{-c}$ as desired. Therefore, the proof of Theorem \ref{thm:density} is complete.
\end{proof}

%%%%%%%%%%%%%%%%%%%%%%%%%%%%%%%%%%%%%%%%%%%%%%%%%%%%%%%%%%%%%%%%%%%%%%%%%%%%%%%%%%%%%%%%%%%%%%%%%%%%%%%%%%%%%%%%%%%%%%%%%%%%%%

\section{Popular common difference}\label{sec:popular}

We plan to prove Theorem \ref{thm:popular} in this section. The guiding philosophy of our approach aligns closely with that of \cite{G-regularity, PPS}. Specifically, given a function $f$ our goal is to find a regular function $g$ and a subset $H$ of $[N]$ such that $f$ and $g$ exhibit similar densities of polynomial configurations with common differences in $H$, while $g$ remains nearly constant under shifts by elements of $H$.

To this end, we introduce the following expression in this section. Suppose that $q\leq M\leq N$ are integers, and $P_1,\dots,P_m\in\Z[y]$ are polynomials. Suppose that $f_0,\dots,f_m:\Z\to\C$ are functions supported on the interval $[N]$. Set
$$
\Lambda^{N,M,q}_{P_1,\dots,P_m}(f_0,f_1,\dots,f_m) =  N^{-1} \sum_{x\in\Z} \E_{y\in [M] \atop q|y} f_0(x) f_1(x+P_1(y)) \cdots f_m(x+P_m(y)).
$$

Given a function $f$, the first step towards our goal is to find regular functions $g_1,\dots,g_m$ and set $H =\set{y\in[M]: q|y}$ such that
$$
\Lambda^{N,M,q}_{P_1,\dots,P_m}(f,\dots,f) \approx \Lambda^{N,M,q}_{P_1,\dots,P_m} (f,g_1,\dots,g_m).
$$

\begin{proposition}[Weak regularity lemma-II]\label{prop:regularity2}
Let $P_1,\dots,P_m\in \Z[y]$ be polynomials of $(C,1)$-coefficients and such that $\deg P_1<\dots<\deg P_m$. Let $d = \deg P_m$ and $d_i = \deg P_i$. Let $\eps \in (0,1)$ be real and $N$ a positive integer.
Then either $N\leq \exp(\exp(O_{C,d}(\eps^{-O_d(1)})))$ or the following statement holds.
If $f_0,\dots,f_m:\Z\to\C$ are $1$-bounded functions supported on the interval $[N]$, then there exist $q,M$ with
$$
N^{1/d} \geq M \geq \exp(\exp(-O_{C,d}(\eps^{-O_d(1)})))N^{1/d}, \ \ q\leq  \exp(\exp(O_{C,d} (\eps^{-O_d(1)}))),
$$
and a local factor chain $(\CB_1,\dots,\CB_m)$ on $[N]$ of resolution $(M^{d_1},\dots,M^{d_m})$ and modulus $q$, such that
$$
\bigabs{\Lambda^{N,\eps M,q}_{P_1,\dots,P_m}(f_0,f_1,\dots,f_m) - \Lambda^{N,\eps M,q}_{P_1,\dots,P_m}(f_0, \Pi_{\CB_1}f_1, \dots, \Pi_{\CB_m}f_m)} \leq \eps. 
$$
\end{proposition}

\begin{remark}
%
%Compared to previous approaches in the literature, our method proceeds directly to a local structural approximation. In particular, it provides an effective framework for obtaining such approximations in cases where the configurations are not easily captured by Gowers uniformity norms. We expect that this perspective may lead to further applications to other problems.

The quantitative bounds for $q$ and $M$ come from the bounds in our inverse theorem (Theorem  \ref{thm:inverse}). Roughly speaking, we will iterative construct a sequence $\{q_j\}_{j \geq 0}$ and take $q = q_j$ for some $j \ll \eps^{-O(1)}$. At each iteration step, we apply the inverse theorem to produce $q_{j+1}$ from $q_j$, with $q_{j+1} = q'q_j^b$ for some $q' \ll \eps^{-O(1)}$ and $b = O(1)$. This implies that $q_j$ is of the shape $\eps^{-O(b^j)}$, leading to the doubly exponential bound for $q$, and hence for $M$ as well.

For the nonlinear Roth pattern studied in \cite[Theorem 6.1]{PPS}, one can obtain (single) exponential bound in the popular difference result, thanks to having a version of Proposition \ref{prop:regularity2} where one can take $q_{j+1} = q'q_j$ in each iteration in the proof. Attaining this improvement in the general setting of Theorem \ref{thm:popular} seems to require an inverse theorem with a precise value of $b$, possibly with $b = d$.
\end{remark}

\begin{proof}
We shall apply the energy increment argument to prove this proposition.
For $j \geq 0$, we will iteratively construct a sequence of local factor chains $(\CB_1^{(j)},\dots,\CB_m^{(j)})$ on $[N]$ of resolution $(M_j^{d_1},\cdots,M_j^{d_m})$ and modulus $q_j$, such that $(\CB_1^{(j+1)},\dots,\CB_m^{(j+1)})$ refines $(\CB_1^{(j)},\dots,\CB_m^{(j)})$ for each $j$, where
\begin{align}\label{eq:Mj-qj}
N^{1/d} \geq M_j \geq \exp(\exp(-O_{C,d}((j+1)\eps^{-O_d(1)})))N^{1/d},\ \ q_j \leq \exp(\exp(O_{C,d} ((j+1)\eps^{-O_d(1)}))).
\end{align}
For $j=0$, let $(\CB_1^{(0)},\cdots,\CB_m^{(0)})$ be a local factor chain of resolution $(M_0^{d_1},\cdots,M_0^{d_m})$ with $M_0=N^{1/d}$ and modulus $q_0=1$. The existence of such a local factor chain easily follows from Lemma \ref{lem:local-factor}. If we have, for some $j \geq 0$,
\begin{equation}\label{eq:reg-terminate}
\bigabs{\Lambda^{N,\eps M_j,q_j}_{P_1,\dots,P_m}(f_0,f_1,\cdots,f_m) - \Lambda^{N,\eps M_j,q_j}_{P_1,\dots,P_m}(f_0, \Pi_{\CB_1^{(j)}}f_1, \cdots, \Pi_{\CB_m^{(j)}}f_m)} \leq \eps, 
\end{equation}
then we terminate the iteration process. We will show that the process must be terminated for some $j \ll_{C,d}\eps^{-O(1)}$, which would conclude the proof by taking $\CB_i = \CB_i^{(j)}$, $q = q_j$, and $M = M_j$.

Now suppose that \eqref{eq:reg-terminate} fails. The left-hand side of \eqref{eq:reg-terminate} can be decomposed into the sum of $m$ terms of the form
$$
\Lambda^{N,\eps M_j,q_j}_{P_1,\dots,P_m}(f_0, g_1,\cdots,g_{i-1}, f_i-\Pi_{\CB_i^{(j)}}f_i, g_{i+1},\cdots,g_m),
$$
where $g_{\ell} = \Pi_{\CB_\ell^{(j)}}f_\ell$ for $\ell < i$ and $g_{\ell} = f_{\ell}$ for $\ell > i$. By the pigeonhole principle, there exists at least one $i\in\set{1,\dots,m}$ such that
\[
\bigabs{\Lambda^{N,\eps M_j,q_j}_{P_1,\dots,P_m}(f_0, g_1,\dots, f_i-\Pi_{\CB_i^{(j)}}f_i, \dots,g_m)}\geq \eps/m.
\]
Set $Q_l(y)=P_l(q_jy)$ for all $1\leq l\leq m$. Recalling Definition \ref{def-lip}, the assumption that $P_1,\dots,P_m$ have $(C,1)$-coefficients implies that $Q_1,\dots,Q_m$ have $(C,q_j^2)$-coefficient. Therefore, recalling the notation \eqref{counting-mn}, we can conclude from the above inequality that
\[
\bigabs{ \Lambda^{N, \eps M_j/q_j}_{Q_1,\dots,Q_m} (f_0, g_1,\dots, f_i-\Pi_{\CB_i^{(j)}}f_i, \dots,g_m)}\geq \eps/m.
\]
By Theorem \ref{thm:inverse}, either we have $\eps M_j/q_j \ll_{C,d} (q_j/\eps)^{O_d(1)}$ in which case the bound $N \leq \exp(\exp(O_{C,d}(\eps^{-O_d(1)})))$ follows and we are done, or else there exist positive integers $q'\ll_{C,d} \eps^{-O(1)}$, $b \ll_d 1$, and a $1$-bounded function $\phi: \Z\rightarrow\C$ which is $O_{C,d}((q_j/\eps)^{O(1)}M_j^{-d_i})$-Lipschitz along $q'q_j^b\cdot\Z$,  such that
\begin{equation}\label{eq:reg-correlation}
\bigabs{\sum_{x \in \Z} (f_i-\Pi_{\CB_i^{(j)}}f_i)(x) \phi(x)} \geq \eta N
\end{equation}
for some $\eta \gg_{C,d}\eps^{O_d(1)}$.
Set $q_{j+1} = q'q_j^b$ so that $q_{j+1}$ satisfies the bound in \eqref{eq:Mj-qj}. We may choose $M_{j+1}$ satisfying the bound in \eqref{eq:Mj-qj} while ensuring that 
$$
q_{j+1}M_{j+1}^{d_i} \leq \frac{M_j^{d_i}}{10}
$$
and $\phi$ is $\eta/(10M_{j+1}^{d_i})$-Lipschitz along $q_{j+1}\cdot\Z$. 
Moreover, for each $1 \leq i \leq m$ we have (with the convention that $d_0=0$)
$$
10 M_j^{d_{i-1}} \leq  10N^{d_{i-1}/d} \leq M_{j+1}^{d_i}
$$
thanks to the lower bound for $M_{j+1}$ in \eqref{eq:Mj-qj}, unless $N \leq \exp(\exp(O_{C,d}(\eps^{-O(1)})))$ in which case we are done. By the two inequalities above, we may apply Lemma \ref{lem:local-factor-chain} to find a local factor chain $(\CB_1^{(j+1)},\cdots,\CB_m^{(j+1)})$ of resolution $(M_{j+1}^{d_1},\cdots,M_{j+1}^{d_m})$ and modulus $q_{j+1}$ which refines $(\CB_1^{(j)},\cdots,\CB_m^{(j)})$.

Since the Lipschitz property of $\phi$ implies that
$$
|\phi(x+q_{j+1}y) - \phi(x)| \leq \frac{1}{5}\eta \text{ for } |y| \leq 2M_{j+1}^{d_i},
$$
the left-hand side of \eqref{eq:reg-correlation} can be written as
$$
\sum_{P} \phi(x_P)\sum_{x \in P} (f_i-\Pi_{\CB_i^{(j)}}f_i)(x) + \sum_{P}\sum_{x \in P}(f_i-\Pi_{\CB_i^{(j)}}f_i)(x) (\phi(x) - \phi(x_P)),
$$
where the outer summation is over all atoms $P$ of $\CB_i^{(j+1)}$ and $x_P$ is an arbitrary element in $P$. In the second term above, each summand is bounded by $\eta/2$ in absolute value since $x-x_P = q_{j+1}y$ for some $|y| \leq 2M_{j+1}^{d_i}$. It follows from inequality (\ref{eq:reg-correlation}) that
$$
\frac{1}{2}\eta N \leq \sum_{P} |\phi(x_P)|\Bigabs{\sum_{x \in P} (f_i-\Pi_{\CB_i^{(j)}}f_i)(x)} \leq \sum_P \Bigabs{\sum_{x \in P} (f_i-\Pi_{\CB_i^{(j)}}f_i)(x)}.
$$
Since $\CB_i^{(j+1)}$ refines $\CB_i^{(j)}$, we have 
$$
\sum_{x \in P} (f_i-\Pi_{\CB_i^{(j)}}f_i)(x) = \sum_{x \in P} \Pi_{\CB_i^{(j+1)}}(f_i-\Pi_{\CB_i^{(j)}}f_i)(x) =  \sum_{x \in P} \Big( \Pi_{\CB_i^{(j+1)}}f_i(x) - \Pi_{\CB_i^{(j)}}f_i(x)\Big),
$$
and hence
$$
\|\Pi_{\CB_i^{(j+1)}}f_i - \Pi_{\CB_i^{(j)}}f_i\|_1 \gg \eta N \gg_{C,d}\eps^{O_d(1)}N.
$$
By orthogonality (Lemma \ref{pythogoras}) and then Cauchy-Schwarz, it follows that
$$
\|\Pi_{\CB_i^{(j+1)}}f_i\|_2 - \|\Pi_{\CB_i^{(j)}}f_i\|_2  =   \|\Pi_{\CB_i^{(j+1)}}f_i - \Pi_{\CB_i^{(j)}}f_i\|_2 \gg_{C,d} \eps^{O_d(1)}N,
$$
and hence
$$
\sum_{\ell=1}^m\|\Pi_{\CB_{\ell}^{(j+1)}}f_\ell\|_2 - \sum_{\ell=1}^m\|\Pi_{\CB_\ell^{(j)}}f_\ell\|_2 \geq \|\Pi_{\CB_i^{(j+1)}}f_i\|_2 - \|\Pi_{\CB_i^{(j)}}f_i\|_2  \gg_{C,d} \eps^{O_d(1)}N.
$$
Hence \eqref{eq:reg-terminate} must hold for some $j  \ll_{C,d} \eps^{-O_d(1)}$. This concludes the proof.
\end{proof}

\begin{proof}[Proof of Theorem \ref{thm:popular}]
Without loss of generality, we may assume that $\deg P_1 < \cdots < \deg P_m$. Let $d = \deg P_m$. Choose a constant $C$ such that $P_1,\cdots,P_m$ have $(C,1)$-coefficients.
By Proposition \ref{prop:regularity2} applied with $f_0=f_1=\cdots=f_m=1_A$, we find integers $q, M$ with
$$
N^{1/d} \geq M \geq \exp(\exp(-O_{C,d}(\eps^{-O_d(1)})))N^{1/d}, \ \ q\leq  \exp(\exp(O_{C,d}(\eps^{-O_d(1)}))),
$$
and a local factor chain $(\CB_1,\cdots,\CB_m)$ on $[N]$ of resolution $(M^{d_1},\cdots,M^{d_m})$ and modulus $q$, such that
$$
\bigabs{\Lambda^{N,\eps M,q}_{P_1,\dots,P_m}(1_A,1_A,\cdots,1_A) - \Lambda^{N,\eps M,q}_{P_1,\dots,P_m}(1_A, \Pi_{\CB_1}1_A, \cdots, \Pi_{\CB_m}1_A)} \leq \eps. 
$$
We will show that
\begin{equation}\label{eq:popular-lower-bound}
\Lambda^{N,\eps M,q}_{P_1,\dots,P_m}(1_A, \Pi_{\CB_1}1_A, \cdots, \Pi_{\CB_m}1_A) \geq \delta^{m+1}-O_{C,d}(\eps).
\end{equation}
This would imply that the expression
$$
 \Lambda_{P_1,\cdots,P_m}^{N,\eps M,q}(1_A,\cdots,1_A) = \frac{1}{N} \sum_{x \in \Z} \E_{\substack{y \leq \eps M \\ q\mid y}} 1_A(x) 1_A(x+P_1(y))\cdots 1_A(x+P_m(y))
$$
is at least $\delta^{m+1} - O_{C,d}(\eps)$.
Hence, by pigeonhole principle,  there exists a positive integer $y \leq \eps M$ with $q \mid y$ such that 
$$
\sum_{x \in \Z} 1_A(x) 1_A(x+P_1(y))\cdots 1_A(x+P_m(y)) \geq (\delta^{m+1}-O_{C,d}(\eps))N,
$$
which concludes the proof.

To prove \eqref{eq:popular-lower-bound}, note that the left-hand side is
$$
\frac{1}{N} \sum_{x \in \Z} \E_{y \leq \eps M \atop q \mid y} 1_A(x) (\Pi_{\CB_1}1_A)(x+P_1(y)) \cdots (\Pi_{\CB_m}1_A)(x+P_m(y)).
$$
For $y \leq \eps M$ and $q\mid y$, we have $q \mid P_i(y)$ and $|P_i(y)| \ll_{C,d}  \eps M^{d_i}$ for each $i$. Thus for such values of $y$, the elements $x, x+P_i(y)$ lie in the same atom of $\CB_i$ for all but $O_{C,d}(\eps N)$ values of $x \in [N]$. It follows that the left-hand side of \eqref{eq:popular-lower-bound} is
$$
\frac{1}{N} \sum_{x \in \Z} 1_A(x) (\Pi_{\CB_1}1_A)(x) \cdots (\Pi_{\CB_m}1_A)(x) + O_{C,d}(\eps).
$$
The desired lower bound follows from the lemma below.
\end{proof}

\begin{lemma}[Combinatorics of projections]\label{combinatorics of projections}
Let $X$ be a finite set of integers and let $f : X \to [0,1]$  be a function with $\E_{x \in X} f(x) = \delta$. Suppose that $\CB_1,\cdots,\CB_m$ are factors of $X$ such that $\CB_i$ refines $\CB_{i+1}$ for each $1 \leq i <m$. Then
$$
\E_{x \in X}  (\Pi_{\CB_1}f)(x) \cdots (\Pi_{\CB_m}f)(x) \geq \delta^{m}.
$$
\end{lemma} 

\begin{proof}

We first claim that for any atom  $P\in \CB_m$ and letting $\delta_P= \E_{x\in P} f(x)$, the following inequality holds
\[
\E_{x\in P}  (\Pi_{\CB_1}f)(x) \cdots (\Pi_{\CB_m}f)(x) \geq \delta_P^m.
\]
Assume this claim, one  can deduce from  H\"older's inequality that
\begin{align*}
\E_{x \in X}  (\Pi_{\CB_1}f)(x) \cdots (\Pi_{\CB_m}f)(x)& = \E_{P\in \CB_m} \E_{x\in P}	(\Pi_{\CB_1}f)(x) \cdots (\Pi_{\CB_m}f)(x)\\
& \geq \E_{P\in \CB_m}\delta_P^m = \E_{P\in \CB_m} (\E_{x\in P} f(x))^m\geq (\E_{x\in X} f(x))^m=\delta^m.
\end{align*}
Thus, it remains only to prove the claim. We proceed by induction on $m$. For the base case $m=1$, it follows immediately from the definition that
\[
\E_{x\in P} \Pi_{\CB_1} f(x) =\E_{x\in P} f(x) =\delta_P.
\]
Now we suppose that $m\geq 2$, and assume the inductive hypothesis holds for $m-1$, that is, if  $Q\in \CB_{m-1}$ and $\delta_Q=\E_Q(f)$ is the density of $f$ on $Q$, then 
\[
\E_{x\in Q} (\Pi_{\CB_1}f)(x) \cdots (\Pi_{\CB_{m-1}}f)(x)\geq \delta_Q^{m-1}.
\]
Since $\Pi_{\CB_m}f$ takes the value $\delta_P$ on the atom $P\in \CB_m$,  we have
\[
\E_{x\in P}  (\Pi_{\CB_1}f)(x) \cdots (\Pi_{\CB_m}f)(x) =\delta_P \E_{x\in P}  (\Pi_{\CB_1}f)(x) \cdots (\Pi_{\CB_{m-1}}f)(x).
\]
Since $\CB_{m-1}$ refines $\CB_m$, we may decompose $P$ as a disjoint union of atoms in $\CB_{m-1}$, i.e. $P=\sqcup_{Q\in \CB_{m-1}\cap P} Q$, Using the inductive hypothesis and  H\"older's inequality, we get 
\begin{align*}
\E_{x\in P}  (\Pi_{\CB_1}f)(x) \cdots (\Pi_{\CB_{m-1}}f)(x) &=\E_{Q\in \CB_{m-1}\cap P}\E_{x\in Q} (\Pi_{\CB_1}f)(x) \cdots (\Pi_{\CB_{m-1}}f)(x)\\
&\geq \E_{Q\in \CB_{m-1}\cap P} \delta_Q^{m-1} \geq (\E_{Q\in \CB_{m-1}\cap P}\E_{x\in Q} f(x))^{m-1}\geq \delta_P^{m-1}.
\end{align*}
Combining the above two  inequalities, we conclude that
\[
\E_{x\in P}  (\Pi_{\CB_1}f)(x) \cdots (\Pi_{\CB_m}f)(x)\geq \delta_P \cdot \delta_P^{m-1} = \delta_P^m,
\]
as claimed.

\end{proof}

%%%%%%%%%%%%%%%%%%%%%%%%%%%%%%%%%%%%%%%%%%%%%%%%%%%%%%%%%%%%%%%%%%%%%%%%%%%%%%%%%%%%%%%%%%%%%%%%%%%%%%%%%%%%%%%%%%%%%%%%%%%%%%

%\renewcommand{\refname}{\normalsize References}
%\bibliographystyle{plain}
%\bibliography{bib}

\bibliographystyle{plain}  % bib style

\bibliography{md_references.bib}

\end{document}